\documentclass{amsart}

\newtheorem{theorem}[equation]{Theorem}
\newtheorem{lemma}[equation]{Lemma}
\newtheorem{corollary}[equation]{Corollary}
\newtheorem{proposition}[equation]{Proposition}

\numberwithin{equation}{section}

\usepackage{amscd}
\usepackage{amsmath}
\usepackage{amssymb}

\begin{document}

\title{On the zeta function of a projective complete intersection}
\author{Alan Adolphson}
\address{Department of Mathematics\\
Oklahoma State University\\
Stillwater, Oklahoma 74078}
\email{adolphs@math.okstate.edu}
\thanks{The first author was supported in part by NSF Grant DMS-0070510} 
\author{Steven Sperber}
\address{School of Mathematics\\
University of Minnesota\\
Minneapolis, Minnesota 55455}
\email{sperber@math.umn.edu}
\date{}
\begin{abstract}
We compute a basis for the $p$-adic Dwork cohomology of a smooth
complete intersection in projective space over a finite field and use
it to give $p$-adic estimates for the action of Frobenius on this
cohomology.  In particular, we prove that the Newton polygon of the
characteristic polynomial of Frobenius lies on or above the associated
Hodge polygon.  This result was first proved by B. Mazur using
crystalline cohomology.  
\end{abstract}
\maketitle

\section{Introduction}

Let $X\subseteq{\bf P}^n$ be a smooth complete intersection of
codimension $r$ over a finite field ${\bf F}_q$ of characteristic $p$.
The zeta function of $X$ may be written in the form 
\[ Z(X/{\bf F}_q,t) =
\frac{P(t)^{(-1)^{n-r-1}}}{(1-t)(1-qt)\cdots(1-q^{n-r}t)}, \]
where $P(t)\in 1+t{\bf Z}[t]$.  The reciprocal roots of $P(t)$ are
units at all primes except the archimedian primes and those primes
lying over $p$.  At any archimedian prime, they have absolute value
$q^{(n-r)/2}$ by Deligne\cite{DE1}; the $p$-adic Newton polygon of
$P(t)$ lies over its Hodge polygon by Mazur\cite{MA1,MA2}.  (The Hodge
polygon is determined by the multidegree of the complete intersection.)

The hypersurface case of Mazur's result was originally proved by
Dwork\cite{DW2}.  Adolphson-Sperber\cite{AS1}, using a toric approach,
extended Dwork's method to (among other cases) the case where $X$ is
the intersection of smooth hypersurfaces of degrees prime to $p$
that meet transversally and which also meet all coordinate
varieties transversally (this latter condition can be ensured by a
coordinate change if $q$ is sufficiently large), but technical
difficulties prevented the 
extension of Dwork's method to general complete intersections until
now.  Recently, we gave an algebraic description of the Jacobian ring
of a complete intersection (see
\cite{AS2}).  This provides the algebraic basis for a proof of Mazur's
theorem by a generalization of Dwork's work.

We describe our results.  Suppose that $X$ is defined by nonconstant
homogeneous polynomials $f_1,\dots,f_r\in{\bf F}_q[x_0,x_1,\dots,x_n]$ of
degrees $d_1,\dots,d_r$, respectively.  Put
\[ F=\sum_{i=1}^r y_if_i\in{\bf F}_q[x_0,\dots,x_n,y_1,\dots,y_r]. \]
We consider the complex $(\Omega^\bullet_{{\bf F}_q[x,y]/{\bf
F}_q},dF\wedge)$, where $dF\in \Omega^1_{{\bf F}_q[x,y]/{\bf F}_q}$
denotes the exterior derivative of $F$ and the boundary map sends
$\omega\in\Omega^k_{{\bf F}_q[x,y]/{\bf F}_q}$ to $dF\wedge\omega\in
\Omega^{k+1}_{{\bf F}_q[x,y]/{\bf F}_q}$.  This complex has a
bigrading defined by taking $(\Omega^k_{{\bf F}_q[x,y]/{\bf
F}_q})^{(e_1,e_2)}$ to be the ${\bf F}_q$-span of those differential forms 
\[ x_0^{a_0}\cdots x_n^{a_n}y_1^{b_1}\cdots y_r^{b_r}\,dx_{i_1}\cdots
dx_{i_l} dy_{j_1}\cdots dy_{j_m} \]
with $l+m=k$,
\[ e_1 = \sum_{i=0}^n a_i + l -\sum_{j=1}^r b_jd_j - d_{j_1}-\cdots-
d_{j_m}, \]
and
\[ e_2 = b_1+\cdots+b_r+m \]
i.~e., we are defining the bidegrees of $x_i$ and $dx_i$ to be $(1,0)$
and the bidegrees of $y_j$ and $dy_j$ to be $(-d_j,1)$.  Since the
boundary map $\omega\mapsto dF\wedge\omega$ is bihomogeneous of
bidegree $(0,1)$, there is an induced bigrading on cohomology that we
denote by
\[ H^k(\Omega^\bullet_{{\bf F}_q[x,y]/{\bf
F}_q},dF\wedge)^{(e_1,e_2)}. \]
Note that as ${\bf F}_q[x,y]$-module, $H^{n+r+1}(\Omega^\bullet_{{\bf
F}_q[x,y]/{\bf F}_q},dF\wedge) = {\bf F}_q[x,y]/I$, where $I$ is the
ideal 
\[ I = (f_1,\dots,f_r,\sum_{j=1}^r y_j\frac{\partial f_j}{\partial
x_0},\dots, \sum_{j=1}^r y_j\frac{\partial f_j}{\partial x_n}). \]
When $f_1=\dots=f_r=0$ defines a smooth complete intersection $X$ in ${\bf
P}^n$, we proved in \cite[Theorem 1.6]{AS2} that
\[ H^{n+r+1}(\Omega^\bullet_{{\bf F}_q[x,y]/{\bf
F}_q},dF\wedge)^{(0,e)} = 0\quad\text{for $e<r$ or $e>n$}. \]
For $r\leq e\leq n$, put
\begin{multline*}
h_e = \dim_{{\bf F}_q} H^{n+r+1}(\Omega^\bullet_{{\bf F}_q[x,y]/{\bf
F}_q},dF\wedge)^{(0,e)} \\
-\begin{cases} 1 & \text{if $p\mid d_1\cdots d_r$, $n+r$ is odd, and
$e=(n+r+1)/2$} \\ 0 & \text{otherwise.} \end{cases}
\end{multline*}
We gave an explicit formula for the polynomial $\sum_{e=r}^n h_et^e$
in \cite[section 5]{AS2}.  In particular, it depends only on $n$, $r$,
and the multidegree $(d_1,\dots,d_r)$ of $X$ (see
\cite[Corollary~1.14]{AS2} for additional properties of the $h_e$).  

By a {\it monomial form\/} we mean an $(n+r+1)$-form of the type
\[ x_0^{a_0}\cdots x_n^{a_n}y_1^{b_1}\cdots y_r^{b_r}\,dx_0\cdots
dx_n\,dy_1\cdots dy_r. \]
In the exceptional case where $p\mid d_1\cdots d_r$, $n+r$ is odd, and
$e=(n+r+1)/2$, let $\xi^{(e)}_l$, $l=1,\dots,h_e$ be monomial forms
such that
\[ \{\xi_l^{(e)}\mid l=1,\dots,h_e\}\cup\{\tau_{(n+r+1)/2}\} \]
is a basis for $H^{n+r+1}(\Omega^\bullet_{{\bf F}_q[x,y]/{\bf
F}_q},dF\wedge)^{(0,e)}$, where $\tau_{(n+r+1)/2}$ is defined in the
proof of Lemma 4.10 below.  In the nonexceptional cases, for
$e=r,\dots,n$, let $\xi^{(e)}_l$, $l=1,\dots,h_e$ be a basis of
monomial forms for 
\[ H^{n+r+1}(\Omega^\bullet_{{\bf F}_q[x,y]/{\bf
F}_q},dF\wedge)^{(0,e)}. \]
Our main result is that the monomial forms
$\{\xi^{(e)}_l\mid e=r,\dots,n,\;l=1,\dots,h_e\}$
are a basis for the Dwork cohomology group
$H^{n+r+1}(\Omega^\bullet_{b,0})$ of $X$. (The definition of Dwork
cohomology is recalled in section 2; a formal statement of this result
is given in Theorem 5.18 below.)  For (generic) smooth affine and toric
complete intersections, Bourgeois\cite{B1,B2} has given an explicit
isomorphism between Dwork cohomology and rigid cohomology.  It seems
likely that his results extend to the projective case, so that the
image of our basis under Bourgeois's isomorphism is a basis for the
rigid cohomology of smooth projective complete intersections.

Let ${\rm ord}_q$ denote the $p$-adic valuation normalized by the
condition ${\rm ord}_q\;q = 1$.  As an application of our work we
prove the following, which is the main result of Mazur(\cite{MA1,MA2}).  
\begin{theorem}
Suppose that $f_1=\dots=f_r=0$ defines a smooth complete intersection $X$
in ${\bf P}^n$.  Then the Newton polygon of $P(t)$ with respect to
${\rm ord}_q$ lies on or above the Newton polygon with respect to
${\rm ord}_q$ of the polynomial  
\[ \prod_{e=0}^{n-r}(1-q^et)^{h_{e+r}}. \]
\end{theorem}

{\it Remark}.  The Newton polygon of
$\prod_{e=0}^{n-r}(1-q^et)^{h_{e+r}}$ is equal to the
middle-dimensional primitive Hodge polygon of $X$ (see
\cite[section 1]{AS2}). 

The outline of the paper is as follows.  In section 2, we review
Dwork's $p$-adic cohomology theory.  In section 3, we define a
$p$-adic filtration on the Dwork complex and use it to compute some of
the cohomology.  We introduce some auxiliary complexes in section 4
and compute their cohomology.  The results of section 4 are applied in
section 5 to finish the computation of the cohomology of the Dwork
complex (Theorem 5.18).  In section 6 we study the action of Frobenius
on cohomology and prove Theorem 1.1.  In section 7, we have collected
as an appendix some general results relating the cohomology of a
$p$-adic complex to the cohomology of its reduction mod $p$.  These
results are used at various points in the paper.

\section{Dwork cohomology}

We retain the notation of the previous section.  Let $X'\subseteq{\bf
A}^{n+1}$ be the affine variety defined by the equations
$f_1=\cdots=f_r=0$.  Let $N_m$ (resp.\ $N'_m$) be the number of ${\bf
F}_{q^m}$-rational points of $X$ (resp.\ $X'$).  Then 
\begin{align*}
Z(X/{\bf F}_q;t)& = \exp\biggl(\sum_{m=1}^\infty
N_m\frac{t^m}{m}\biggr) \\
Z(X'/{\bf F}_q;t)& = \exp\biggl(\sum_{m=1}^\infty
N'_m\frac{t^m}{m}\biggr) 
\end{align*}
The obvious relation
\[ N_m=\frac{N_m'-1}{q^m-1} \]
is equivalent to
\begin{equation}
Z(X'/{\bf F}_q;t)=\frac{Z(X/{\bf F}_q;qt)}{(1-t)Z(X/{\bf F}_q;t)}.
\end{equation}
If we write
\[ Z(X/{\bf F}_q;t)=\frac{P(t)^{(-1)^{n-r-1}}}{(1-t)(1-qt)\cdots
(1-q^{n-r}t)}, \]
then (2.1) implies
\begin{equation}
Z(X'/{\bf
F}_q;t)=\frac{1}{1-q^{n-r+1}t}\biggl(\frac{P(t)}{P(qt)}\biggr)^{(-1)^{n-r}}.
\end{equation}

This zeta function is closely related to the $L$-function of a certain
exponential sum.  Let $\Psi:{\bf F}_q\rightarrow{\bf C}^{\times}$ be a
nontrivial additive character and consider the exponential sums
\[ S_m=\sum_{(x,y)\in{\bf A}^{n+1+r}({\bf F}_{q^m})} \Psi({\rm
Tr}_{{\bf F}_{q^m}/{\bf F}_q}(F(x,y))). \]
The associated $L$-function is
\[ L({\bf A}^{n+1+r},\Psi,F;t)=\exp\biggl(\sum_{m=1}^\infty
S_m\frac{t^m}{m}\biggr). \]
It is easily seen that
\[ \sum_{y\in{\bf A}^r({\bf F}_q)} \Psi({\rm
Tr}_{{\bf F}_{q^m}/{\bf F}_q}(F(x,y))) = \begin{cases}
q^{mr}& \text{if $f_1(x)=\cdots=f_r(x)=0$}, \\
0& \text{otherwise}, 
\end{cases} \]
hence
\[ L({\bf A}^{n+1+r},\Psi,F;t)=Z(X'/{\bf F}_q;q^rt). \]
By (2.2), we thus have
\begin{equation}
L({\bf A}^{n+1+r},\Psi,F;t)^{(-1)^{n+r}} = (1-q^{n+1}t)^{(-1)^{n+r-1}}
\frac{P(q^rt)}{P(q^{r+1}t)}.
\end{equation}

We recall the construction of the Dwork complex whose cohomology
describes this $L$-function.  Let ${\bf Q}_p$ be the field of $p$-adic numbers,
$\zeta_p$ a primitive $p$-th root of unity, and $\Lambda_1={\bf
  Q}_p(\zeta_p)$.  The field $\Lambda_1$ is a totally ramified
extension of ${\bf Q}_p$ of degree $p-1$.  Write $q=p^a$ and let $K$
be the unramified extension of ${\bf Q}_p$ of degree~$a$.  Set
$\Lambda_0=K(\zeta_p)$.  The Frobenius automorphism $x\mapsto x^p$ of
${\rm Gal}({\bf F}_q/{\bf F}_p)$ lifts to a generator $\tau$ of ${\rm
  Gal}(\Lambda_0/\Lambda_1)(\simeq {\rm Gal}(K/{\bf Q}_p))$ by requiring
$\tau(\zeta_p)=\zeta_p$.  Let $\Lambda$ be the completion of an
algebraic closure of $\Lambda_0$.  Denote by ``ord'' the additive
valuation on $\Lambda$ normalized by ${\rm ord}\;p=1$ and by ``${\rm
  ord}_q$'' the additive valuation normalized by ${\rm ord}_q\;q=1$.

Let $E(t)$ be the Artin-Hasse exponential series:
\[ E(t)=\exp\biggl(\sum_{i=0}^{\infty} \frac{t^{p^i}}{p^i}\biggr). \]
Let $\gamma\in\Lambda_1$ be a solution of $\sum_{i=0}^{\infty}
t^{p^i}/p^i=0$ satisfying ${\rm ord}\;\gamma=1/(p-1)$ and consider
\begin{equation}
\theta(t)=E(\gamma t)=\sum_{i=0}^{\infty}
\lambda_it^i\in\Lambda_1[[t]].
\end{equation}
The series $\theta(t)$ is a splitting function in Dwork's
terminology\cite[section 4a]{DW1}.  Its coefficients satisfy
\begin{equation}
{\rm ord}\;\lambda_i\geq i/(p-1).
\end{equation}

We consider the following spaces of $p$-adic functions.  Let $b$ be a
positive rational number and choose a positive integer $M$ such that
both $Mb/(p(p-1))$ and $M/(p-1)$ are integers.  Let $\pi$ be such that
\begin{equation}
\pi^{M}=p
\end{equation}
and put $\tilde{\Lambda}_1=\Lambda_1(\pi)$,
$\tilde{\Lambda}_0=\Lambda_0(\pi)$.  The element $\pi$ is a uniformizing
parameter for the rings of integers of $\tilde{\Lambda}_1$ and
$\tilde{\Lambda}_0$.  We extend $\tau\in {\rm Gal}(\Lambda_0/\Lambda_1)$
to a generator of ${\rm Gal}(\tilde{\Lambda}_0/\tilde{\Lambda}_1)$ by
requiring $\tau(\pi)=\pi$.  For $v=(v_1,\ldots,v_r)\in{\bf R}^r$, we
put $|v|=v_1+\cdots+v_r$.  Define
\begin{multline}
C(b)= \\
\biggl\{ \sum_{(u,v)\in{\bf N}^{n+1+r}} A_{u,v}\pi^{Mb|v|}x^uy^v \mid
  \text{$A_{u,v}\in\tilde{\Lambda}_0$ and $A_{u,v}\rightarrow 0$ as
  $(u,v)\rightarrow\infty$} \biggr\}. 
\end{multline}

We construct a de Rham-type complex using $C(b)$.  Let
\[ \Omega_b^k=\bigoplus_{l+m=k}\bigoplus_{\substack{0\leq
i_1<\cdots<i_l\leq n\\1\leq j_1<\cdots<j_m\leq r}} C(b)\,
dx_{i_1}\cdots dx_{i_l} dy_{j_1}\cdots dy_{j_m}. \]
Define an exterior derivative $d:\Omega^k_b\to\Omega^{k+1}_b$ by
linearity and the formula
\begin{multline*} 
d(\xi\,dx_{i_1}\cdots dx_{i_l}\, dy_{j_1}\cdots dy_{j_m}) = \\
\biggl(\sum_{i=0}^n \frac{\partial\xi}{\partial x_i}\,dx_i +
\sum_{j=1}^r \frac{\partial\xi}{\partial y_j}\,dy_j\biggr)\wedge
dx_{i_1}\cdots dx_{i_l}\, dy_{j_1}\cdots dy_{j_m} 
\end{multline*}
for $\xi\in C(b)$.  Let $\hat{f}_i = \sum_u \hat{a}_{u,i}x^u\in K[x]$
be the Teichm\"uller lifting of $f_i$, let
\begin{equation}
\hat{F}=\sum_{i=1}^r y_i\hat{f}_i(x) = \sum_{u,i} \hat{a}_{u,i}x^uy_i
\in K[x,y]
\end{equation}
be the Teichm\"{u}ller lifting of $F$, and put
$\gamma_l=\sum_{i=0}^l \gamma^{p^i}/p^i$.  From the definition of
$\gamma$ we have
\begin{equation}
{\rm ord}\;\gamma_l\geq \frac{p^{l+1}}{p-1}-l-1.
\end{equation}
Set
\begin{equation}
H = \sum_{l=0}^{\infty} \gamma_l\hat{F}^{\tau^l}(x^{p^l},y^{p^l}) =
\sum_{l=0}^\infty \gamma_l\sum_{i=1}^r
y_i^{p^l}\hat{f}_i^{\tau^l}(x^{p^l}). 
\end{equation}
It follows from (2.9) that $\partial H/\partial x_i,\partial
H/\partial y_j\in C(b)$, so $dH\in\Omega^1_b$.  Define
$D:\Omega^k_b\to\Omega^{k+1}_b$ by
\[ D(\omega) = \pi^{Mb}\gamma^{-1}(d\omega + dH\wedge\omega). \]
(The reason for the normalizing factor $\pi^{Mb}\gamma^{-1}$ will
become apparent in section 3.)  We thus obtain a complex
$(\Omega^\bullet_b,D)$.  

We define the Frobenius operator on this complex.  Set (see (2.8))
\begin{align}
G(x)& = \prod_{u,i} \theta(\hat{a}_{u,i}x^uy_i), \\
G_0(x)& = \prod_{j=0}^{a-1}\prod_{u,i} \theta((\hat{a}_{u,i}x^uy_i)^{p^j}).
\end{align}
The estimate (2.5) implies that $G\in C(b)$ for all $b<1/(p-1)$ and
$G_0\in C(b)$ for all $b<p/q(p-1)$.  Define an operator $\psi$ on
formal power series by 
\begin{equation}
\psi\biggl(\sum_{(u,v)\in{\bf
N}^{n+1+r}}A_{u,v}x^uy^v\biggr)=\sum_{(u,v)\in{\bf N}^{n+1+r}} 
A_{pu,pv}x^uy^v.
\end{equation}
It is clear that $\psi(C(b))\subseteq C(pb)$.  For $0<b<p/(p-1)$,
let $\alpha=\psi^a\circ G_0$ be the composition 
\[ C(b)\hookrightarrow C(b/q)\xrightarrow{G_0}
C(b/q)\xrightarrow{\psi^a} C(b), \] 
where the middle arrow is multiplication by $G_0$.  Then $\alpha$ is a
completely continuous $\tilde{\Lambda}_0$-linear endomorphism of
$C(b)$.   We shall also consider
$\beta=\tau^{-1}\circ\psi\circ G$, which is a completely continuous
$\tilde{\Lambda}_1$-linear (or $\tilde{\Lambda}_0$-semilinear)
endomorphism of~$C(b)$.  Note that $\alpha=\beta^a$.  

Define a map $\alpha_\bullet:\Omega^\bullet_b\to\Omega^\bullet_b$ by
additivity and the formula ($k=l+m$)
\begin{multline}
\alpha_k(\xi\,dx_{i_1}\cdots dx_{i_l}\,dy_{j_1}\cdots dy_{j_m}) = \\
\frac{q^{n+1+r-k}}{x_{i_1}\cdots x_{i_l}y_{j_1}\cdots
y_{j_m}}\alpha(x_{i_1}\cdots x_{i_l}y_{j_1}\cdots y_{j_m}\xi)\,
dx_{i_1}\cdots dx_{i_l}\,dy_{j_1}\cdots dy_{j_m}. 
\end{multline}
Then $\alpha_\bullet$ is a map of complexes (see, e.~g.,
\cite[equation (2.11)]{AS3}).  Similarly, we define a map of complexes
(see \cite[equation (2.12)]{AS3})
$\beta_\bullet:\Omega^\bullet_b\to\Omega^\bullet_b$ by additivity and the
formula
\begin{multline}
\beta_k(\xi\,dx_{i_1}\cdots dx_{i_l}\,dy_{j_1}\cdots dy_{j_m}) = \\
\frac{p^{n+1+r-k}}{x_{i_1}\cdots x_{i_l}y_{j_1}\cdots
y_{j_m}}\beta(x_{i_1}\cdots x_{i_l}y_{j_1}\cdots y_{j_m}\xi)\,
dx_{i_1}\cdots dx_{i_l}\,dy_{j_1}\cdots dy_{j_m}. 
\end{multline}
The Dwork trace formula, as formulated by Robba\cite{RO}, then gives
\begin{equation}
L({\bf A}^{n+1+r},\Psi,F;t)=\prod_{k=0}^{n+r+1} \det(I-t\alpha_k\mid
\Omega_b^k)^{(-1)^{k+1}}.
\end{equation}

So far we have followed closely the description of Dwork's theory as
given, for example, in \cite{AS3}.  However, since the $f_i$ are
homogeneous polynomials, we can replace $\Omega_b^\bullet$ by a
smaller complex that is more easily analyzed.  For $s\in{\bf Z}$, let
$\Omega_{b,s}^k$ be the subspace of $\Omega_b^k$ spanned by those
$k$-forms ($l+m=k$)
\[ \sum_{u,v}A_{u,v}\pi^{Mb|v|}x^uy^v\,dx_{i_1}\cdots
dx_{i_l}\,dy_{j_1}\cdots dy_{j_m} \] 
with 
\begin{equation}
\sum_{i=0}^n u_i +l-\sum_{j=1}^r v_jd_j - d_{j_1}-\cdots-d_{j_m} =s
\end{equation}
for all $(u,v)$ with $A_{u,v}\neq 0$.
There is an obvious decomposition as $\tilde{\Lambda}_0$-vector space
\[ \Omega_b^k=\bigoplus_{s\in{\bf Z}} \Omega_{b,s}^k. \]
Note that every monomial
$x^uy^v$ that appears in any of $F,\hat{F},G,G_0$ with nonzero
coefficient satisfies $\sum_{i=0}^n u_i=\sum_{j=1}^r v_jd_j$.  
It then follows easily that
\begin{equation}
D(\Omega_{b,s}^k)\subseteq \Omega_{b,s}^{k+1}
\end{equation}
and that
\begin{equation} 
\alpha_k(\Omega_{b,s}^k)\subseteq \Omega_{b,s/q}^k.
\end{equation}
In particular, if $q\nmid s$, $\alpha_k$ is zero on $\Omega^k_{b,s}$;
so for fixed $s\neq 0$, $\alpha_k^N$ is zero on $\Omega^k_{b,s}$ for all
$N\gg 0$.  We thus see that $\alpha_\bullet$ is stable on the subcomplex
$(\Omega_{b,0}^\bullet,D)$ of $(\Omega_b^\bullet,D)$  while 
\[ \det(I-t\alpha_k\mid \Omega_b^k/\Omega_{b,0}^k)=1. \]
Equation (2.16) and Serre\cite[Proposition 9]{SE} now imply
\[ L({\bf A}^{n+1+r},\Psi,F;t)=\prod_{k=0}^{n+r+1}
\det(I-t\alpha_k\mid \Omega_{b,0}^k)^{(-1)^{k+1}}, \]
and passing to cohomology gives
\begin{equation}
L({\bf A}^{n+1+r},\Psi,F;t)=\prod_{k=0}^{n+r+1}
\det(I-t\alpha_k\mid H^k(\Omega_{b,0}^\bullet,D))^{(-1)^{k+1}}.
\end{equation}

\section{The $\pi$-adic filtration on $(\Omega^\bullet_b,D)$}

The $p$-adic Banach space $C(b)$ has a decreasing filtration
$\{F^sC(b)\}_{s=-\infty}^{\infty}$ defined by setting
\[ F^sC(b)=\bigg\{\sum_{(u,v)\in{\bf N}^{n+1+r}} A_{u,v}\pi^{Mb|v|}x^uy^v\in
C(b) \mid A_{u,v}\in\pi^s{\mathcal O}_{\tilde{\Lambda}_0} \text{ for all
  $u,v$}\bigg\}, \]
where ${\mathcal O}_{\tilde{\Lambda}_0}$ denotes the ring of integers
of $\tilde{\Lambda}_0$.  We extend this to a filtration on
$\Omega^{\bullet}_b$ by defining
\[ F^s\Omega^k_b=\bigoplus_{l+m=k}\bigoplus_{\substack{0\leq
  i_1<\cdots<i_l\leq n \\  1\leq j_1<\cdots<j_m\leq r}} \pi^{Mbm}F^sC(b)\,
  dx_{i_1}\cdots dx_{i_l}\, dy_{j_1}\cdots  dy_{j_m}. \]
We assume from now on that $1/(p-1)<b<p/(p-1)$.  A calculation shows that
under this condition, $D(F^s\Omega^k_b)\subset F^s\Omega^{k+1}_b$
(this is the reason for introducing the normalizing factor
$\pi^{Mb}\gamma^{-1}$ in section 2), hence $(\Omega^\bullet_b,D)$ is a
filtered complex.  We shall use this filtration to compute the
cohomology of $(\Omega^\bullet_{b,0},D)$ and then use (2.20) and (2.3)
to estimate the Newton polygon of $P(t)$.   

Consider the map $F^0\Omega^k_b\to\Omega^k_{{\bf F}_q[x,y]/{\bf F}_q}$
defined by additivity and the formula
\begin{multline*} \sum_{(u,v)\in{\bf N}^{n+1+r}}
A_{u,v}\pi^{Mb(|v|+m)}x^uy^v\,dx_{i_1}\cdots dx_{i_l}\,dy_{j_1}\cdots
dy_{j_m} \mapsto \\
\sum_{(u,v)\in{\bf N}^{n+1+r}}
\bar{A}_{u,v}x^uy^v\,dx_{i_1}\cdots dx_{i_l}\,dy_{j_1}\cdots dy_{j_m},
\end{multline*}
where $\bar{A}_{u,v}$ denotes the reduction of $A_{u,v}$ modulo the
maximal ideal of $\tilde{\Lambda}_0$. (Since $A_{u,v}\to 0$ as
$(u,v)\to\infty$, the sum on the right-hand side is finite.)  It is
clear that this map is surjective with kernel $F^1\Omega^k_b$, hence
we get an isomorphism 
\begin{equation}
F^0\Omega^k_b/F^1\Omega^k_b \cong \Omega^k_{{\bf F}_q[x,y]/{\bf F}_q}.
\end{equation}
For $\omega\in F^0\Omega^k_b$ we have
\begin{align*}
D(\omega) &\equiv \pi^{Mb}\gamma^{-1}dH\wedge\omega
\pmod{F^1\Omega^{k+1}_b} \\
&\equiv \pi^{Mb}\biggl( \sum_{i=0}^n\frac{\partial\hat{F}}{\partial
x_i}\,dx_i + \sum_{j=1}^r \frac{\partial\hat{F}}{\partial
y_j}\,dy_j\biggr)\wedge\omega \pmod{F^1\Omega^{k+1}_b} \\
&\equiv \pi^{Mb}d\hat{F}\wedge\omega \pmod{F^1\Omega^{k+1}_b}.
\end{align*}
It follows that under the isomorphism (3.1), $D$ is identified with
$dF\wedge$, i.~e., there is an isomorphism of complexes
\[ (F^0\Omega^\bullet_b/F^1\Omega^\bullet_b,D)\cong
(\Omega^\bullet_{{\bf F}_q[x,y]/{\bf F}_q},dF\wedge). \]
Since multiplication by $\pi^s$ defines an isomorphism of complexes
\begin{equation}
(F^0\Omega^\bullet_b,D)\cong (F^s\Omega^\bullet_b,D), 
\end{equation}
we have in fact isomorphisms for all $s\in{\bf Z}$
\begin{equation}
(F^s\Omega^\bullet_b/F^{s+1}\Omega^\bullet_b,D)\cong
(\Omega^\bullet_{{\bf F}_q[x,y]/{\bf F}_q},dF\wedge).
\end{equation}

Referring to the bigrading defined in the introduction, we set
\[ (\Omega^k_{{\bf F}_q[x,y]/{\bf F}_q})^{(0)} = \bigcup_{e=0}^\infty 
(\Omega^k_{{\bf F}_q[x,y]/{\bf F}_q})^{(0,e)}\subseteq \Omega^k_{{\bf
F}_q[x,y]/{\bf F}_q}. \] 
Since $dF\in (\Omega^1_{{\bf F}_q[x,y]/{\bf
F}_q})^{(0)}$, it follows that $((\Omega^\bullet_{{\bf F}_q[x,y]/{\bf
F}_q})^{(0)},dF\wedge)$ is a subcomplex of $(\Omega^\bullet_{{\bf
F}_q[x,y]/{\bf F}_q},dF\wedge)$.  The isomorphism (3.1) induces
\[ F^0\Omega^k_{b,0}/F^1\Omega^k_{b,0}\cong (\Omega^k_{{\bf
F}_q[x,y]/{\bf F}_q})^{(0)}, \]
so as above we get isomorphisms for all $s$
\begin{equation}
(F^s\Omega^\bullet_{b,0}/F^{s+1}\Omega^\bullet_{b,0},D)\cong
((\Omega^\bullet_{{\bf F}_q[x,y]/{\bf F}_q})^{(0)},dF\wedge).
\end{equation}

Our general approach will be as follows.  The results of \cite{AS2}
describe the cohomology of $((\Omega^\bullet_{{\bf F}_q[x,y]/{\bf
F}_q})^{(0)},dF\wedge)$.  We shall use the isomorphism (3.4) to infer a
description of the cohomology of $(F^s\Omega^\bullet_{b,0},D)$, and from
that we shall calculate the cohomology of $(\Omega^\bullet_{b,0},D)$.
This explicit description will lead to the desired $p$-adic estimates
for the action of Frobenius on the cohomology of $\Omega^\bullet_{b,0}$.

For example, by \cite[Theorem~1.6]{AS2} we have 
\begin{equation}
H^k((\Omega^\bullet_{{\bf F}_q[x,y]/{\bf F}_q})^{(0)},dF\wedge) = 0
\quad\text{for $k\neq 2r,n+r,n+r+1$.}
\end{equation}
An application of Proposition 7.1 then gives the
following result. 
\begin{proposition}
If the equations $f_1=\dots=f_r=0$ define a smooth complete
intersection $X$ in ${\bf P}^n$, then for all $s\in{\bf Z}$
\[ H^k(F^s\Omega^\bullet_{b,0},D) = 0\quad\text{for $k\neq
2r,n+r,n+r+1$.} \] 
\end{proposition}

We next compute $H^{2r}(F^s\Omega^\bullet_{b,0},D)$.  Define
\begin{align}
\xi_1 &= \pi^{Mb}\gamma^{-1}\sum_{j=1}^r
d_1\cdots\hat{d}_j\cdots d_r\biggl(\sum_{l=0}^\infty \frac{\gamma_l}{p^l}
d(\hat{f}_j^{\tau^l}(x^{p^l})) \wedge d(y_j^{p^l})\biggr) \\ \nonumber
&= \pi^{Mb}\gamma^{-1}\sum_{j=1}^r
d_1\cdots\hat{d}_j\cdots d_r\biggl(\sum_{l=0}^\infty \gamma_l y_j^{p^l-1}
d(\hat{f}_j^{\tau^l}(x^{p^l}))\,dy_j\biggr)\in\Omega^2_{b,0}.
\end{align}
For $1<k\leq r$ we define inductively
\begin{align}
\xi_k &= \frac{(-1)^{k-1}}{kd_1\cdots d_r}\xi_1\wedge\xi_{k-1} \\ \nonumber
&= \frac{(-1)^{k(k-1)/2}}{k!\,(d_1\cdots
d_r)^{k-1}}\xi_1\wedge\cdots\wedge\xi_1 \quad\text{($k$ factors).}
\end{align}
Put
\[ \sigma_j = \sum_{l=0}^\infty
\gamma_l y_j^{p^l-1}d(\hat{f}_j^{\tau^l}(x^{p^l})). \]
Then we have explicitly
\begin{equation}
\xi_k = (\pi^{Mb}\gamma^{-1})^k \sum_{1\leq j_1<\dots<j_k\leq r}
\biggl(\prod_{j\not\in\{j_1,\dots,j_k\}} d_j\biggr)\,
\sigma_{j_1}\wedge\cdots\wedge \sigma_{j_k}\wedge
dy_{j_1}\wedge\cdots\wedge dy_{j_k}. 
\end{equation}

If $\xi\in\Omega^i_{b,0}$, then $d(\xi\wedge\eta) = d\xi\wedge\eta +
(-1)^i\xi\wedge d\eta$.  Since $d\circ d = 0$, it follows from (3.7)
and (3.8) that $d\xi_k = 0$ for all $k$.  From (2.10) we have
\[ dH = \sum_{j=1}^r\biggl(y_j\sigma_j + \biggl(\sum_{l=0}^\infty
\gamma_lp^ly_j^{p^l-1}
\hat{f}_j^{\tau^l}(x^{p^l})\biggr)dy_j\biggr). \]
It now follows immediately from (3.9) that $dH\wedge\xi_r = 0$, hence
\begin{equation}
D(\xi_r) = 0.
\end{equation}
It is straightforward to check that $\xi_k\in F^0\Omega^{2k}_{b,0}$
and that under the isomorphism (3.4) (with $s=0$), $\xi_k$ is
mapped to
\[ \bar{\xi}_k = \sum_{1\leq j_1<\dots<j_k\leq
r}\biggl(\prod_{j\not\in\{j_1,\dots,j_k\}} d_j\biggr)\,
df_{j_1}\wedge\cdots\wedge df_{j_k}\wedge dy_{j_1}\wedge\cdots\wedge
dy_{j_k}. \] 
By \cite[Theorem 1.6]{AS2}, if $r<n$ then
\begin{equation}
H^{2r}((\Omega^\bullet_{{\bf F}_q[x,y]/{\bf F}_q})^{(0)},dF\wedge) =
{\bf F}_q\cdot[\bar{\xi}_r].
\end{equation}
\begin{proposition}
If $r<n$ and the equations $f_1=\dots=f_r=0$ define a smooth complete
intersection $X$ in ${\bf P}^n$, then
$H^{2r}(F^s\Omega^\bullet_{b,0},D)$ is a free ${\mathcal
O}_{\tilde{\Lambda}_0}$-module with basis $[\pi^s\xi_r]$. 
\end{proposition}

{\it Proof}.  The isomorphism (3.4) is equivalent to the exactness of
\begin{equation}
0\to (F^s\Omega^\bullet_{b,0},D)\xrightarrow{\pi}
(F^s\Omega^\bullet_{b,0},D)\to ((\Omega^\bullet_{{\bf F}_q[x,y]/{\bf
F}_q})^{(0)},dF\wedge)\to 0,
\end{equation}
where the second arrow is multiplication by $\pi$.  By (3.5), the
associated long exact sequence gives an exact sequence 
\begin{multline*}
0\to H^{2r}(F^s\Omega^\bullet_{b,0},D)\xrightarrow{\pi}
H^{2r}(F^s\Omega^\bullet_{b,0},D)\to H^{2r}((\Omega^\bullet_{{\bf
F}_q[x,y]/{\bf F}_q})^{(0)},dF\wedge) \\
\xrightarrow{\delta}H^{2r+1}(F^s\Omega^\bullet_{b,0},D)\xrightarrow{\pi}
H^{2r+1}(F^s\Omega^\bullet_{b,0},D).
\end{multline*}
Since $\bar{\xi}_r\in(\Omega^{2r}_{{\bf F}_q[x,y]/{\bf F}_q})^{(0)}$
is the image of $\pi^s\xi_r\in F^s\Omega^{2r}_{b,0}$ under the
isomorphism (3.4), we have $\delta([\bar{\xi}_r]) = [D(\pi^s\xi_r)] = 0$
by (3.10).  By (3.11), the connecting homomorphism $\delta$ is the
zero map.  It follows that multiplication by $\pi$ is injective on
$H^{2r}(F^s\Omega^\bullet_{b,0},D)$ and
$H^{2r+1}(F^s\Omega^\bullet_{b,0},D)$.  The assertion of the
proposition now follows immediately from Proposition 7.2.

Propositions 3.6 and 3.12 compute $H^k(F^s\Omega^\bullet_{b,0},D)$ for
$0\leq k<n+r$.  To get information about this cohomology when
$k=n+r,n+r+1$, we need to introduce some related complexes.  

\section{The $\theta$-map}

Define $\theta:\Omega^k_b\to\Omega^{k-1}_b$ by $C(b)$-linearity and
the formula
\begin{multline}
\theta(dx_{i_1}\dots dx_{i_l}\,dy_{j_1}\dots dy_{j_m}) \\
= \sum_{s=1}^l (-1)^{s-1} x_{i_s}\,dx_{i_1}\cdots\widehat{dx}_{i_s}\cdots
dx_{i_l}\,dy_{j_1}\cdots dy_{j_m} \\
+\sum_{t=1}^m (-1)^{l+t-1}(-d_{j_t}y_{j_t})\,dx_{i_1}\cdots
dx_{i_l}\,dy_{j_1}\cdots \widehat{dy}_{j_t}\cdots dy_{j_m}.
\end{multline}
One has $\theta^2=0$, $\theta(\Omega^k_{b,s})\subseteq
\Omega^{k-1}_{b,s}$, and $\theta(F^s\Omega^k_b)\subseteq
F^s\Omega^{k-1}_b$.
\begin{proposition}
The sequence
\[ 0\to F^s\Omega^{n+r+1}_{b,0}\xrightarrow{\theta}
F^s\Omega^{n+r}_{b,0}\xrightarrow{\theta}\cdots\xrightarrow{\theta} 
F^s\Omega^0_{b,0}\to\pi^s{\mathcal O}_{\tilde{\Lambda}_0}\to 0 \]
is exact for all $s\in{\bf Z}$ (where the map
$F^s\Omega^0_{b,0}\to\pi^s{\mathcal O}_{\tilde{\Lambda}_0}$ sends a
power series to its constant term).
\end{proposition}

{\it Proof}.  Define $\theta:(\Omega^k_{{\bf F}_q[x,y]/{\bf
F}_q})^{(0)} \to (\Omega^{k-1}_{{\bf F}_q[x,y]/{\bf F}_q})^{(0)}$ by
${\bf F}_q[x,y]$-linearity and the formula (4.1).  As in (3.4), one
has isomorphisms of complexes
\begin{equation}
(F^s\Omega^\bullet_{b,0}/F^{s+1}\Omega^\bullet_{b,0},\theta)\cong
((\Omega^\bullet_{{\bf F}_q[x,y]/{\bf F}_q})^{(0)},\theta) 
\end{equation}
for all $s$.  By \cite[Proposition 4.6]{AS2} we have 
\[ H_k((\Omega^\bullet_{{\bf F}_q[x,y]/{\bf F}_q})^{(0)},\theta) =
\begin{cases} 0 & \text{if $k>0$} \\ {\bf F}_q & \text{if $k=0$.}
\end{cases} \]
It now follows from Propositions 7.1 and 7.2 that
\[ H_k(F^s\Omega^\bullet_{b,0},\theta) = \begin{cases} 0 & \text{if
$k>0$} \\ \pi^s{\mathcal O}_{\tilde{\Lambda}_0} & \text{if $k=0$,}
\end{cases} \]
which is the assertion of the proposition.

For $k\geq 0$, put 
\begin{equation}
F^s\widetilde{\Omega}_{b,0}^k = \theta(F^s\Omega^{k+1}_{b,0})
\quad(=\text{ker($\theta:F^s\Omega^k_{b,0}\to F^s\Omega^{k-1}_{b,0})$ if
$k>0$.})
\end{equation}
It is straightforward to check that
\begin{equation}
\theta\circ D + D\circ\theta = 0
\end{equation}
on $\Omega^\bullet_{b,0}$, hence $D(F^s\widetilde{\Omega}^k_{b,0})\subset
F^s\widetilde{\Omega}^{k+1}_{b,0}$.  We thus obtain a complex
$(F^s\widetilde{\Omega}^\bullet_{b,0},D)$.  

We define a related complex $F^s\widehat{\Omega}^\bullet_{b,0}$ as
follows.  Let
\[ F^s\widehat{\Omega}^0_{b,0} = \pi^s{\mathcal O}_{\tilde{\Lambda}_0}
\] 
and let $F^s\widehat{\Omega}^k_{b,0} = F^s\widetilde{\Omega}^{k-1}_{b,0}$
for $k\geq 1$.  We define the boundary map
$F^s\widehat{\Omega}^k_{b,0}\to F^s\widehat{\Omega}^{k+1}_{b,0}$ to be $0$
if $k=0$ and $-D$ if $k\geq 1$.  Thus 
\begin{equation}
H^0(F^s\widehat{\Omega}^\bullet_{b,0}) = \pi^s{\mathcal
O}_{\tilde{\Lambda}_0} 
\end{equation}
and 
\begin{equation}
H^k(F^s\widehat{\Omega}^\bullet_{b,0}) =
H^{k-1}(F^s\widetilde{\Omega}^\bullet_{b,0})\quad\text{for $k\geq 1$.}
\end{equation}

Define maps $F^s\Omega^k_{b,0}\to F^s\widehat{\Omega}^k_{b,0}$ as follows.
For $k=0$, take the map $F^s\Omega^0_{b,0}\to\pi^s{\mathcal
O}_{\tilde{\Lambda}_0}$ that sends a power series to its constant
term.  By Proposition 4.2, this map defines an isomorphism
$F^s\Omega^0_{b,0}/\theta(F^s\Omega^1_{b,0})\cong \pi^s{\mathcal 
O}_{\tilde{\Lambda}_0}$.  For $k\geq 1$, take the map
$\theta:F^s\Omega^k_{b,0}\to F^s\widehat{\Omega}^k_{b,0}$.  It follows from
Proposition 4.2 that we have a short exact sequence of complexes
\begin{equation}
0\to F^s\widetilde{\Omega}^\bullet_{b,0}\to F^s\Omega^\bullet_{b,0}\to
F^s\widehat{\Omega}^\bullet_{b,0}\to 0,
\end{equation}
and this gives rise to an exact cohomology sequence
\begin{equation}
\cdots\to H^k(F^s\widetilde{\Omega}^\bullet_{b,0})\to
H^k(F^s\Omega^\bullet_{b,0})\to H^k(F^s\widehat{\Omega}^\bullet_{b,0})
\xrightarrow{\delta}
H^{k+1}(F^s\widetilde{\Omega}^\bullet_{b,0})\to\cdots.
\end{equation}

Our goal in this section is to describe the cohomology groups
$H^k(F^s\widetilde{\Omega}^\bullet_{b,0})$
($=H^{k+1}(F^s\widehat{\Omega}^\bullet_{b,0})$ by (4.7)) for $k<n+r$.
We begin by describing some differential forms $\eta_k\in
F^0\widehat{\Omega}^{2k}_{b,0}$
($=F^0\widetilde{\Omega}^{2k-1}_{b,0}$ for $k>0$) that will play a key
role in what follows.  

\begin{lemma}
Set $\eta_0 = 1\in F^0\widehat{\Omega}^0_{b,0}$.  For $k\geq 1$, there
exist differential forms $\eta_k\in
F^0\widetilde{\Omega}^{2k-1}_{b,0}$ such that \\
{\bf (a)} $[\eta_k]=\delta([\eta_{k-1}])$, where $\delta$ is the
connecting homomorphism in $(4.9)$, and \\
{\bf (b)} $\bar{\eta}_k\in (\Omega^{2k-1}_{{\bf F}_q[x,y]/{\bf
F}_q})^{(0,k)}$, where $\bar{\eta}_k$ is the image of $\eta_k$
under the isomorphism~$(3.4)$.
\end{lemma}

{\it Proof}.  The proof is by induction on $k$.  Define $\eta_1 =
D(1)$.  Then $\bar{\eta}_1 = dF\in 
(\Omega^1_{{\bf F}_q[x,y]/{\bf F}_q})^{(0,1)}$, so the result is true
for $k=1$.  Suppose the result true for some $k\geq 1$.  Now
$\theta(\eta_k) = 0$, so $\theta(\bar{\eta}_k) = 0$.  It follows from
\cite[Proposition 4.6]{AS2} that there exists $\zeta_k\in 
(\Omega^{2k}_{{\bf F}_q[x,y]/{\bf F}_q})^{(0,k)}$ such that
$\theta(\zeta_k) = \bar{\eta}_k$.  By Proposition 7.1
we can choose $\tau_k\in F^0\Omega^{2k}_{b,0}$ such that
$\theta(\tau_k) = \eta_k$ and $\bar{\tau}_k = \zeta_k$.  Define
\begin{equation}
\eta_{k+1} = D(\tau_k).
\end{equation}
Then $\delta([\eta_k]) = \eta_{k+1}$ by the definition of the
connecting homomorphism.  Furthermore, by (4.11), 
\[ \bar{\eta}_{k+1}=dF\wedge\zeta_k\in (\Omega^{2k+1}_{{\bf
F}_q[x,y]/{\bf F}_q})^{(0,k+1)}, \]
and by induction the proof is complete.

\begin{proposition}
For $0\leq k<2r$, if $k$ is even, then 
\[ H^k(F^s\widetilde{\Omega}^\bullet_{b,0}) =
H^{k+1}(F^s\widehat{\Omega}^\bullet_{b,0}) = 0 \]
and if $k$ is odd, then
\[ H^k(F^s\widetilde{\Omega}^\bullet_{b,0}) =
H^{k+1}(F^s\widehat{\Omega}^\bullet_{b,0}) \cong {\mathcal
O}_{\tilde{\Lambda}_0} \cdot [\pi^s\eta_{(k+1)/2}] \]
(a free, rank-one ${\mathcal O}_{\tilde{\Lambda}_0}$-module with basis
$[\pi^s\eta_{(k+1)/2}]$).  
\end{proposition}

{\it Proof}.  By Proposition 3.6, $H^k(F^s\Omega^\bullet_{b,0}) = 0$
for $k<2r$, so 
\begin{equation}
H^0(F^s\widetilde{\Omega}^\bullet_{b,0})) = 0
\end{equation}
and the connecting homomorphism $\delta$ in (4.9) gives isomorphisms
\begin{equation}
H^k(F^s\widehat{\Omega}^\bullet_{b,0})\cong
H^{k+1}(F^s\widetilde{\Omega}^\bullet_{b,0}) \quad\text{for $0\leq k<2r-1$.}
\end{equation}
Equation (4.13) establishes the proposition for $k=0$, and (4.6) and
(4.14) establish it for $k=1$.  Using (4.7), we may regard (4.14) as
isomorphisms 
\begin{equation}
H^k(F^s\widetilde{\Omega}^\bullet_{b,0})\cong
H^{k+2}(F^s\widetilde{\Omega}^\bullet_{b,0}) \quad\text{for $0\leq
k<2r-2$.} 
\end{equation}
Using (4.15) and Lemma 4.10, the proposition follows by induction from
the cases $k=0$ and $k=1$.  

\begin{lemma}
For $k=2,\dots,r$, $\theta(\xi_k) = (-1)^{k-1}D(\xi_{k-1})$.
\end{lemma}

{\it Proof}.  A straightforward calculation shows that
\begin{equation}
\theta(\xi_1) = \pi^{Mb}\gamma^{-1}(d_1\cdots d_r)dH = D(d_1\cdots d_r).
\end{equation}
It is convenient to set $\xi_0 = d_1\dots d_r$, so that (4.17) establishes
the case $k=1$ of the lemma and we can proceed by induction on $k$.
Let $k\geq 2$ and suppose the assertion of the lemma is true for
$1,\dots,k-1$.  If $\omega_1$ is an $l$-form, then
\[ \theta(\omega_1\wedge\omega_2) = \theta(\omega_1)\wedge\omega_2 +
(-1)^l\omega_1\wedge\theta(\omega_2), \]
so by (3.8) and induction we have
\begin{align}
\theta(\xi_k) &= \frac{(-1)^{k-1}}{kd_1\cdots d_r}\biggl(
\theta(\xi_1)\wedge\xi_{k-1} + \xi_1\wedge\theta(\xi_{k-1})\biggr) \\
 &= \frac{(-1)^{k-1}}{kd_1\cdots d_r}\biggl(
 D(\xi_0)\wedge\xi_{k-1} +\xi_1\wedge(-1)^{k-2}D(\xi_{k-2})\biggr). \nonumber
\end{align}
We observed earlier that $d\xi_l=0$, hence $D(\xi_l) = (\pi^{Mb}\gamma^{-1})
dH\wedge\xi_l$ for all $l$.  Substituting this into (4.18) and using
(3.8) gives 
\begin{align*}
\theta(\xi_k) &= \frac{(-1)^{k-1}}{kd_1\cdots
d_r}\biggl( \pi^{Mb}\gamma^{-1}(kd_1\cdots d_r)dH\wedge\xi_{k-1} \biggr) \\
 &= (-1)^{k-1}D(\xi_{k-1}),
\end{align*}
which proves the lemma.

The following result is key to describing the cohomology groups
$H^k(F^s\widetilde{\Omega}^\bullet_{b,0})$ for $k\geq 2r$.  
\begin{proposition}
Let $r<n$.  Relative to the bases $[\pi^s\xi_r]$ for
$H^{2r}(F^s\Omega^\bullet_{b,0})$ and $[\pi^s\eta_r]$ for
$H^{2r}(F^s\widehat{\Omega}^\bullet_{b,0})$, the map
\[ \theta:H^{2r}(F^s\Omega^\bullet_{b,0})\to
H^{2r}(F^s\widehat{\Omega}^\bullet_{b,0}) \]
is multiplication by $(-1)^{r(r-1)/2}d_1\cdots d_r$.
\end{proposition}

{\it Proof}.  It suffices to prove the assertion when $s=0$.  We prove
inductively that for $k=1,\dots,r$,
\begin{equation}
\theta(\xi_k) = (-1)^{k(k-1)/2}(d_1\cdots d_r)\eta_k +
D(\theta(\sigma_k))
\end{equation}
for some $\sigma_k\in F^0\Omega^{2k-1}_{b,0}$.  The assertion of the
proposition then follows by taking $k=r$ in (4.20).  By (4.17),
\[ \theta(\xi_1) = (d_1\cdots d_r)D(1) = (d_1\cdots d_r)\eta_1, \]
so the assertion is true for $k=1$.  Assume inductively that (4.20) holds
for some~$k$, $1\leq k<r$.  Choose $\tau_k\in F^0\Omega^{2k}_{b,0}$,
as in the proof of Lemma 4.10, so that $\theta(\tau_k) =
\eta_k$.  Substitution into (4.20) then gives (since $D\circ\theta =
-\theta\circ D$)
\[ \theta(\xi_k - (-1)^{k(k-1)/2}(d_1\cdots d_r)\tau_k + D(\sigma_k))
= 0. \]
By Proposition 4.2, there exists $\sigma_{k+1}\in
F^0\Omega^{2k+1}_{b,0}$ such that
\begin{equation}
\xi_k = (-1)^{k(k-1)/2}(d_1\cdots d_r)\tau_k - D(\sigma_k) +
\theta((-1)^k\sigma_{k+1}).
\end{equation}
From Lemma 4.16 we have
\begin{equation}
\theta(\xi_{k+1}) = (-1)^kD(\xi_k).
\end{equation}
Substituting (4.21) into (4.22) now gives
\[ \theta(\xi_{k+1}) = (-1)^{k(k+1)/2}(d_1\cdots d_r)D(\tau_k) +
D(\theta(\sigma_{k+1})). \]
Since $\eta_{k+1} = D(\tau_k)$ (see (4.11)), this is just (4.20) with $k$
replaced by $k+1$, so by induction the proof is complete.

\begin{proposition}
Let $2r\leq k<n+r$. \\
{\bf (a)} If $(p,d_1\cdots d_r) = 1$, then 
\[ H^k(F^s\widetilde{\Omega}^\bullet_{b,0}) =
H^{k+1}(F^s\widehat{\Omega}^\bullet_{b,0}) = 0. \]
{\bf (b)} if $p\mid d_1\cdots d_r$, then
\[ H^k(F^s\widetilde{\Omega}^\bullet_{b,0}) =
H^{k+1}(F^s\widehat{\Omega}^\bullet_{b,0}) \cong \begin{cases} {\mathcal
O}_{\tilde{\Lambda}_0}/(d_1\cdots d_r){\mathcal O}_{\tilde{\Lambda}_0}
& \text{if $k$ is odd,} \\ 0 & \text{if $k$ is even.} \end{cases} \]
\end{proposition}

{\it Proof:}  If $r=n$ there is nothing to prove, so suppose $r<n$.
From (4.9) and Proposition 4.12 we have an exact sequence 
\begin{equation}
0\to H^{2r}(F^s\widetilde{\Omega}^\bullet_{b,0})\to
H^{2r}(F^s\Omega^\bullet_{b,0})\xrightarrow{\theta}
H^{2r}(F^s\widehat{\Omega}^\bullet_{b,0}).
\end{equation}
By Proposition 4.19, $\theta$ is injective, so
\begin{equation}
H^{2r}(F^s\widetilde{\Omega}^\bullet_{b,0})=
H^{2r+1}(F^s\widehat{\Omega}^\bullet_{b,0})=0 .
\end{equation}
If $r=n-1$ there is nothing more to prove, so suppose $r<n-1$.  By
Proposition~3.6 we have $H^{2r+1}(F^s\Omega^\bullet_{b,0})=0$, so
(4.9) gives an exact sequence
\begin{equation}
H^{2r}(F^s\Omega^\bullet_{b,0})\xrightarrow{\theta}
H^{2r}(F^s\widehat{\Omega}^\bullet_{b,0}) \xrightarrow{\delta}
H^{2r+1}(F^s\widetilde{\Omega}^\bullet_{b,0})\to 0. 
\end{equation}
Proposition 4.19 then implies
\begin{equation}
H^{2r+1}(F^s\widetilde{\Omega}^\bullet_{b,0})=
H^{2r+2}(F^s\widehat{\Omega}^\bullet_{b,0}) \cong {\mathcal 
O}_{\tilde{\Lambda}_0} / (d_1\cdots d_r){\mathcal
O}_{\tilde{\Lambda}_0}.
\end{equation}

Assume first that $(p,d_1\cdots d_r)=1$, so that by (4.25) and (4.27)
we have
\[ H^{2r}(F^s\widetilde{\Omega}^\bullet_{b,0}) =
H^{2r+1}(F^s\widetilde{\Omega}^\bullet_{b,0}) = 0. \]
Suppose that for some $k$, $2r\leq k<n+r-2$, we have proved
\[ H^k(F^s\widetilde{\Omega}^\bullet_{b,0}) =
H^{k+1}(F^s\widetilde{\Omega}^\bullet_{b,0}) = 0. \]
Then by (4.7) we also have
\[ H^{k+1}(F^s\widehat{\Omega}^\bullet_{b,0}) =
H^{k+2}(F^s\widehat{\Omega}^\bullet_{b,0}) = 0. \]
The exact sequence (4.9) then implies
\[ H^{k+2}(F^s\widetilde{\Omega}^\bullet_{b,0}) \cong
H^{k+2}(F^s\Omega^\bullet_{b,0}). \]
But $H^{k+2}(F^s\Omega^\bullet_{b,0})=0$ by Proposition 3.6, so
\[ H^{k+2}(F^s\widetilde{\Omega}^\bullet_{b,0}) =
H^{k+3}(F^s\widehat{\Omega}^\bullet_{b,0}) = 0. \]
Part (a) of the proposition now follows by induction on $k$.

Now assume that $p\mid d_1\cdots d_r$.  Consider first the case of
even $k$.  The assertion
\begin{equation}
H^k(F^s\widetilde{\Omega}^\bullet_{b,0}) =
H^{k+1}(F^s\widehat{\Omega}^\bullet_{b,0}) = 0
\end{equation}
holds for $k=2r$ by (4.25).  Suppose it holds for some even $k$,
$2r\leq k\leq n+r-3$.  By Proposition 3.6 we have 
\[ H^{k+1}(F^s\Omega^\bullet_{b,0}) = H^{k+2}(F^s\Omega^\bullet_{b,0})
= 0, \]
so (4.9) gives an isomorphism
\begin{equation}
H^{k+1}(F^s\widehat{\Omega}^\bullet_{b,0})\cong
H^{k+2}(F^s\widetilde{\Omega}^\bullet_{b,0}). 
\end{equation}
Equation (4.28) then implies that
\begin{equation}
H^{k+2}(F^s\widetilde{\Omega}^\bullet_{b,0}) =
H^{k+3}(F^s\widehat{\Omega}^\bullet_{b,0}) = 0,
\end{equation}
so part (b) of the proposition follows for even $k$ by induction on
$k$.  

Now consider the case of odd $k$.  The assertion
\begin{equation}
H^k(F^s\widetilde{\Omega}^\bullet_{b,0}) =
H^{k+1}(F^s\widehat{\Omega}^\bullet_{b,0}) \cong {\mathcal
O}_{\tilde{\Lambda}_0} / (d_1\cdots d_r) {\mathcal
O}_{\tilde{\Lambda}_0} 
\end{equation}
holds for $k=2r+1$ by (4.27).  Suppose it holds for some odd $k$,
$2r+1\leq k\leq n+r-3$.  The isomorphism (4.29) then gives
\begin{equation}
H^{k+2}(F^s\widetilde{\Omega}^\bullet_{b,0}) =
H^{k+3}(F^s\widehat{\Omega}^\bullet_{b,0}) \cong {\mathcal
O}_{\tilde{\Lambda}_0} / (d_1\cdots d_r) {\mathcal
O}_{\tilde{\Lambda}_0},
\end{equation}
so part (b) of the proposition follows for odd $k$ by induction on $k$
also. 

We describe generators for the torsion modules of Proposition 4.23.  
\begin{proposition}
Suppose $p\mid d_1\cdots d_r$.  For $k$ odd, $2r+1\leq k<n+r$, the
cohomology class $[\pi^s\eta_{(k+1)/2}]$ generates the torsion
module
\[ H^k(F^s\widetilde{\Omega}^\bullet_{b,0}) =
H^{k+1}(F^s\widehat{\Omega}^\bullet_{b,0}) \cong {\mathcal 
O}_{\tilde{\Lambda}_0} / (d_1\cdots d_r) {\mathcal
O}_{\tilde{\Lambda}_0}. \]
\end{proposition}

{\it Proof}.  It suffices to prove the assertion for $s=0$.  The
proof is by induction on $k$.  Consider first $k=2r+1$.  The exact
sequence (4.26) shows that $\delta([\eta_r])$ is a generator
for $H^{2r+1}(F^0\widetilde{\Omega}^\bullet_{b,0})$.  But
$\delta([\eta_r]) = [\eta_{r+1}]$ by Lemma 4.10, so the result is true
for $k=2r+1$.
Suppose inductively that for some odd $k$, $2r+1\leq k<n+r$, the
proposition is true.  
By (4.29), $\delta([\eta_{(k+1)/2}])$ generates
$H^{k+2}(F^0\widetilde{\Omega}^\bullet_{b,0})$.  But by Lemma~4.10,
$\delta([\eta_{(k+1)/2}]) = [\eta_{(k+3)/2}]$, so by induction on $k$
the proof is complete.

\section{Computation of $H^{n+r}(F^s\Omega^\bullet_{b,0})$ and
$H^{n+r+1}(F^s\Omega^\bullet_{b,0})$} 

\begin{lemma}
$H^{n+r}(F^s\Omega^\bullet_{b,0})$ is a free, finitely-generated
${\mathcal O}_{\tilde{\Lambda}_0}$-module.
\end{lemma}

{\it Proof}:  The isomorphism (3.4) is equivalent to the exactness of
the sequence
\begin{equation}
0\to (F^0\Omega^\bullet_{b,0},D)\xrightarrow{\pi}
(F^0\Omega^\bullet_{b,0},D)\to  ((\Omega^\bullet_{{\bf F}_q[x,y]/{\bf
F}_q})^{(0)},dF\wedge)\to 0, 
\end{equation}
where the second arrow is multiplication by $\pi$.  From the
associated sequence of cohomology groups we get the exact sequence
\begin{equation}
H^{n+r-1}((\Omega^\bullet_{{\bf F}_q[x,y]/{\bf F}_q})^{(0)})
\xrightarrow{\delta} H^{n+r}(F^0\Omega^\bullet_{b,0})
\xrightarrow{\pi} H^{n+r}(F^0\Omega^\bullet_{b,0}). 
\end{equation}
If $r\neq n-1$ (so that $n+r-1\neq 2r$), equation (3.5) implies that
multiplication by $\pi$ is injective on
$H^{n+r}(F^0\Omega^\bullet_{b,0})$.  If $r=n-1$, then (5.3) becomes 
\begin{equation}
H^{2r}((\Omega^\bullet_{{\bf F}_q[x,y]/{\bf F}_q})^{(0)})
\xrightarrow{\delta} H^{n+r}(F^0\Omega^\bullet_{b,0})
\xrightarrow{\pi} H^{n+r}(F^0\Omega^\bullet_{b,0}). 
\end{equation}
But by (3.11), this connecting homomorphism $\delta$ is the zero map:
\[ \delta([\bar{\xi}_r]) = [D(\xi_r)] = 0 \]
by (3.10).  Thus in all cases multiplication by $\pi$ is injective on
$H^{n+r}(F^0\Omega^\bullet_{b,0})$, hence
$H^{n+r}(F^0\Omega^\bullet_{b,0})$ is torsion-free.  To show it is
free, we are thus reduced to showing that it is finitely generated.

Since $F^0\widetilde{\Omega}^k_{b,0} = 0$ for $k>n+r$, we get the from
(4.9) the exact sequence
\[ H^{n+r}(F^0\widetilde{\Omega}^\bullet_{b,0})\to
H^{n+r}(F^0\Omega^\bullet_{b,0})\to
H^{n+r}(F^0\widehat{\Omega}^\bullet_{b,0})\to 0. \]
By Proposition 4.23, $H^{n+r}(F^0\widehat{\Omega}^\bullet_{b,0})$ is
finitely generated, so it suffices to show that
$H^{n+r}(F^0\widetilde{\Omega}^\bullet_{b,0})$ is finitely generated.
By (4.7) we have 
\[ H^{n+r}(F^0\widetilde{\Omega}^\bullet_{b,0}) =
H^{n+r+1}(F^0\widehat{\Omega}^\bullet_{b,0}), \]
and from (4.9) we get $H^{n+r+1}(F^0\Omega^\bullet_{b,0})\cong
H^{n+r+1}(F^0\widehat{\Omega}^\bullet_{b,0})$, so we are finally
reduced to showing that $H^{n+r+1}(F^0\Omega^\bullet_{b,0})$ is
finitely generated.  But by \cite[Theorem 1.6]{AS2}
$H^{n+r+1}((\Omega^\bullet_{{\bf F}_q[x,y]/{\bf F}_q})^{(0)})$ is
finitely generated.  We can choose a ``monomial'' basis for it
consisting of $(n+r+1)$-forms of the type
\[ x_0^{a_0}\cdots x_n^{a_n}y_1^{b_1}\cdots y_r^{b_r}\,dx_0\cdots dx_n
dy_1\cdots dy_r. \]
After multiplying by an appropriate power of $\pi$, such a form may
be regarded as an element of $F^0\Omega^{n+r+1}_{b,0}$.  These
normalized forms span $H^{n+r+1}(F^0\Omega^\bullet_{b,0})$ by
Proposition 7.3.  This completes the proof of the
lemma.

\begin{proposition}
If $(p,d_1\cdots d_r) = 1$ or if $n+r$ is
even, then $H^{n+r+1}(F^s\Omega^\bullet_{b,0})$ is a free ${\mathcal
O}_{\tilde{\Lambda}_0}$-module.  
\end{proposition}

{\it Proof}.  It suffices to prove the assertion for $s=0$.  From
(4.9) we get the exact sequence 
\begin{equation}
H^{n+r-1}(F^0\Omega^\bullet_{b,0})\xrightarrow{\theta}
H^{n+r-1}(F^0\widehat{\Omega}^\bullet_{b,0})\to
H^{n+r}(F^0\widetilde{\Omega}^\bullet_{b,0})\to
H^{n+r}(F^0\Omega^\bullet_{b,0}).
\end{equation}
If $r=n$, then $H^{n+r-1}(F^0\widehat{\Omega}^\bullet_{b,0})=0$ by
Proposition 4.12.  If $r=n-1$ and $(p,d_1\cdots d_r) = 1$, then the
first arrow in (5.6) is surjective by Proposition 4.19.  The case
$r=n-1$ cannot occur when $n+r$ is even.  Finally, if $r<n-1$, then 
$H^{n+r-1}(F^0\widehat{\Omega}^\bullet_{b,0})=0$ by Proposition 4.23.
It follows that in all cases, the last arrow in (5.6) is 
injective, so $H^{n+r}(F^0\widetilde{\Omega}^\bullet_{b,0})$ is free
by Lemma 5.1.  But $H^{n+r}(F^0\widetilde{\Omega}^\bullet_{b,0}) =
H^{n+r+1}(F^0\widehat{\Omega}^\bullet_{b,0})$ by (4.7), and from (4.9),
$H^{n+r+1}(F^s\Omega^\bullet_{b,0})\cong
H^{n+r+1}(F^0\widehat{\Omega}^\bullet_{b,0})$.  Thus
$H^{n+r+1}(F^s\Omega^\bullet_{b,0})$ is also free.

\begin{proposition}
If $p\mid d_1\cdots d_r$ and $n+r$ is odd, then the torsion submodule
of $H^{n+r+1}(F^s\Omega^\bullet_{b,0})$ is isomorphic to ${\mathcal
O}_{\tilde{\Lambda}_0}/(d_1\cdots d_r){\mathcal
O}_{\tilde{\Lambda}_0}$ and is generated by the cohomology class
$[\pi^s\tau_{(n+r+1)/2}]$ (where $\tau_{(n+r+1)/2}$ is defined in the proof
of Lemma~$4.10$).  Furthermore, there is an exact sequence 
\[ 0\to \langle[\pi^s\tau_{(n+r+1)/2}]\rangle \to
H^{n+r+1}(F^s\Omega^\bullet_{b,0}) \xrightarrow{\theta}
H^{n+r}(F^s\Omega^\bullet_{b,0})\to 0, \] 
where $\langle[\pi^s\tau_{(n+r+1)/2}]\rangle$ denotes the submodule
generated by $[\pi^s\tau_{(n+r+1)/2}]$.
\end{proposition}

{\it Proof}.  Since $n+r$ is odd, the case $r=n$ is impossible so we
have $r<n$.  From (4.9) and Proposition 4.23(b) we get the exact sequence
\begin{multline}
H^{n+r-1}(F^0\Omega^\bullet_{b,0})\xrightarrow{\theta}
H^{n+r-1}(F^0\widehat{\Omega}^\bullet_{b,0})\xrightarrow{\delta}
H^{n+r}(F^0\widetilde{\Omega}^\bullet_{b,0}) \\
\to H^{n+r}(F^0\Omega^\bullet_{b,0})\to 0.
\end{multline}
By Proposition 4.12 if $r=n-1$ and by Proposition 4.33 if $r<n-1$, it
follows that $H^{n+r-1}(F^0\widehat{\Omega}^\bullet_{b,0})$ is
generated by $[\eta_{(n+r-1)/2}]$.  Hence by Lemma 4.10, the image of
$\delta$ is generated by $[\eta_{(n+r+1)/2}]$.  By Proposition~4.19 if
$r=n-1$ and by Propositions~3.6 and~4.23(b) if $r<n-1$, it follows
that the cokernel of $\theta$ in (5.8) is isomorphic to ${\mathcal
O}_{\tilde{\Lambda}_0}/(d_1\cdots d_r){\mathcal
O}_{\tilde{\Lambda}_0}$.  Thus (5.8) gives an exact sequence
\begin{equation}
0\to\langle[\eta_{(n+r+1)/2}]\rangle\to
H^{n+r}(F^0\widetilde{\Omega}^\bullet_{b,0})\to
H^{n+r}(F^0\Omega^\bullet_{b,0})\to 0
\end{equation}
where $\langle[\eta_{(n+r+1)/2}]\rangle\cong {\mathcal
O}_{\tilde{\Lambda}_0}/(d_1\cdots d_r){\mathcal
O}_{\tilde{\Lambda}_0}$.  By (4.7) we have
$H^{n+r}(F^0\widetilde{\Omega}^\bullet_{b,0}) =
H^{n+r+1}(F^0\widehat{\Omega}^\bullet_{b,0})$, and (4.9) shows that the
map $\theta$ gives an isomorphism 
\[ H^{n+r+1}(F^0\Omega^\bullet_{b,0})
\cong H^{n+r+1}(F^0\widehat{\Omega}^\bullet_{b,0}). \]
Under this isomorphism, the cohomology class $[\tau_{(n+r+1)/2}]\in 
H^{n+r+1}(F^0\Omega^\bullet_{b,0})$ corresponds to
$[\eta_{(n+r+1)/2}]\in H^{n+r+1}(F^0\widehat{\Omega}^\bullet_{b,0})$.
With these identifications, the sequence (5.9) becomes
\[ 0\to \langle[\tau_{(n+r+1)/2}]\rangle\to
H^{n+r+1}(F^0\Omega^\bullet_{b,0}) \xrightarrow{\theta}
H^{n+r}(F^0\Omega^\bullet_{b,0})\to 0, \] 
with $\langle[\tau_{(n+r+1)/2}]\rangle\cong {\mathcal
O}_{\tilde{\Lambda}_0}/(d_1\cdots d_r){\mathcal
O}_{\tilde{\Lambda}_0}$.  This proves the proposition.

{\it Remark}.  The exact sequence (5.2) implies that there is an
isomorphism
\[ H^{n+r+1}(F^0\Omega^\bullet_{b,0})/\pi
H^{n+r+1}(F^0\Omega^\bullet_{b,0}) \cong
H^{n+r+1}((\Omega^\bullet_{{\bf F}_q[x,y]/{\bf F}_q})^{(0)}). \]
From Proposition 5.7, it then follows that the cohomology class
$[\bar{\tau}_{(n+r+1)/2}]$ is not zero in
$H^{n+r+1}((\Omega^\bullet_{{\bf F}_q[x,y]/{\bf F}_q})^{(0)})$.  

In section 1, we defined $h_e$ to be the dimension
of $H^{n+r+1}(\Omega^\bullet_{{\bf F}_q[x,y]/{\bf F}_q})^{(0,e)}$ for
$r\leq e\leq n$  
except in the exceptional case $p\mid d_1\cdots d_r$, $n+r$ odd,
$e=(n+r+1)/2$, where this dimension is $h_e+1$.  Choose elements of
$(\Omega^{n+r+1}_{{\bf F}_q[x,y]/{\bf F}_q})^{(0,e)}$ of the form
\begin{equation}
\bar{\xi}_l^{(e)} = \prod_{i=0}^n x_i^{a_i(e;l)} \prod_{j=1}^r
y_j^{b_j(e;l)}\, dx_0\cdots dx_n dy_1\cdots dy_r
\end{equation}
for $l=1,\dots,h_e$ so that $\{[\bar{\xi}_l^{(e)}]\}_{l=1}^{h_e}$ is a
basis for $H^{n+r+1}(\Omega^\bullet_{{\bf F}_q[x,y]/{\bf
F}_q})^{(0,e)}$ in the nonexceptional cases and
$\{[\bar{\xi}_l^{((n+r+1)/2)}]\}_{l=1}^{h_{(n+r+1)/2}}\cup
\{[\bar{\tau}_{(n+r+1)/2}]\}$ is a basis in the exceptional case.
As noted
in section 1,
\[ H^{n+r+1}((\Omega^\bullet_{{\bf F}_q[x,y]/{\bf F}_q})^{(0)}) =
\bigcup_{e=r}^n H^{n+r+1}(\Omega^\bullet_{{\bf F}_q[x,y]/{\bf
F}_q})^{(0,e)}, \]
so $\{[\bar{\xi}_l^{(e)}]\mid e=r,\dots,n,\;l=1,\dots,h_e\}$ is a
basis for $H^{n+r+1}((\Omega^\bullet_{{\bf F}_q[x,y]/{\bf
F}_q})^{(0)})$ if either $(p,d_1\cdots d_r) = 1$ or $n+r$ is
even, and $\{[\bar{\xi}_l^{(e)}]\mid e=r,\dots,n,\;l=1,\dots,h_e\}\cup
\{[\bar{\tau}_{(n+r+1)/2}]\}$ is a basis if $p\mid d_1\cdots d_r$ and
$n+r$ is odd.  Put 
\begin{equation}
\xi_l^{(e)} = \pi^{Mbe}\prod_{i=0}^n x_i^{a_i(e;l)} \prod_{j=1}^r
y_j^{b_j(e;l)}\, dx_0\cdots dx_n dy_1\cdots dy_r\in F^0\Omega^{n+r+1}_{b,0}.
\end{equation}
The image of $\xi_l^{(e)}$ under the isomorphism (3.4) is
$\bar{\xi}_l^{(e)}$, so by Propositions 5.5 and~7.2 we get the
following result. 
\begin{corollary} 
If $(p,d_1\cdots d_r)=1$ or if $n+r$ is
even, then the cohomology classes $\{[\pi^s\xi_l^{(e)}]\mid
e=r,\dots,n,\;l=1,\dots,h_e\}$ form a basis for
$H^{n+r+1}(F^s\Omega^\bullet_{b,0})$.
\end{corollary}

The following result describes the situation in the exceptional case.
\begin{corollary}
If $p\mid d_1\cdots d_r$ and $n+r$ is odd, then 
\[ H^{n+r+1}(F^s\Omega^\bullet_{b,0}) = H_s\oplus\langle
[\pi^s\tau_{(n+r+1)/2}] \rangle, \]
where $H_s$ is a free submodule of
$H^{n+r+1}(F^s\Omega^\bullet_{b,0})$ with basis the cohomology classes
$\{[\pi^s\xi_l^{(e)}]\mid e=r,\dots,n,\;l=1,\dots,h_e\}$.
\end{corollary}

{\it Proof}.  The fact that $H^{n+r+1}(F^s\Omega^\bullet_{b,0})$ is
the sum of these two submodules follows from Proposition 7.3.  Suppose
there are $c_{e,l}\in {\mathcal O}_{\tilde{\Lambda}_0}$ such that 
\begin{equation}
\sum_{e,l} c_{e,l}[\xi_l^{(e)}] = 0 
\end{equation}
in $H^{n+r+1}(F^0\Omega^\bullet_{b,0})$.  Then in
$H^{n+r+1}((\Omega^\bullet_{{\bf F}_q[x,y]/{\bf F}_q})^{(0)})$ we have
\begin{equation}
\sum_{e,l} \bar{c}_{e,l}[\bar{\xi}_l^{(e)}] = 0.
\end{equation}
If some $c_{e,l}$ is a unit in ${\mathcal O}_{\tilde{\Lambda}_0}$, then
the relation (5.15) is nontrivial, contradicting the definition
of the $\bar{\xi}_l^{(e)}$.  If all $c_{e,l}$ are divisible by $\pi$ but
some $c_{e,l}$ is nonzero, we
can choose $w$ so that $c_{e,l}=\pi^w c'_{e,l}$ for all $e,l$, where
$c'_{e,l}\in {\mathcal O}_{\tilde{\Lambda}_0}$ and some $c'_{e,l}$ is
a unit.  Then (5.14) says that $\sum_{e,l} c'_{e,l}[\xi_l^{(e)}]$ is a
torsion element of $H^{n+r+1}(F^0\Omega^\bullet_{b,0})$, so by Proposition
5.7 we have 
\begin{equation}
\sum_{e,l} c'_{e,l}[\xi_l^{(e)}] = c[\tau_{(n+r+1)/2}]
\end{equation}
for some $c\in {\mathcal O}_{\tilde{\Lambda}_0}$.  But this implies
the nontrivial relation 
\begin{equation}
\sum_{e,l} \bar{c}'_{e,l}[\bar{\xi}_l^{(e)}] =
\bar{c}[\bar{\tau}_{(n+r+1)/2}]
\end{equation}
in $H^{n+r+1}((\Omega^\bullet_{{\bf F}_q[x,y]/{\bf F}_q})^{(0)})$,
again contradicting the definition of the $\bar{\xi}_l^{(e)}$. 

The following result is the key to estimating the Newton polygon.  It
is an immediate consequence of Corollaries 5.12 and 5.13 and
Proposition 7.7.
\begin{theorem}
In all cases, the cohomology classes $\{[\xi_l^{(e)}]\mid
e=r,\dots,n,\;l=1,\dots,h_e\}$ form a basis for
$H^{n+r+1}(\Omega^\bullet_{b,0})$.
\end{theorem}

\section{Frobenius action and Newton polygon}

It follows from Proposition 4.2 that the sequence
\begin{equation}
0\to\Omega^{n+r+1}_{b,0}\xrightarrow{\theta}\Omega^{n+r}_{b,0}
\xrightarrow{\theta}\dots\xrightarrow{\theta}
\Omega^0_{b,0}\to\tilde{\Lambda}_0\to 0
\end{equation}
is exact.  Put
\[ \widetilde{\Omega}^k_{b,0} = \theta(\Omega^{k+1}_{b,0}) =
\bigcup_{s\in{\bf Z}} F^s\widetilde{\Omega}^k_{b,0} \]
and define
\begin{align*}
\widehat{\Omega}^0_{b,0} &= \tilde{\Lambda}_0 \\
\widehat{\Omega}^k_{b,0} &= \widetilde{\Omega}^{k-1}_{b,0} =
\bigcup_{s\in{\bf Z}} F^s\widehat{\Omega}^k_{b,0} \quad\text{for
$k\geq 1$.}
\end{align*}
The boundary maps $F^s\widetilde{\Omega}^k_{b,0}\to
F^s\widetilde{\Omega}^{k+1}_{b,0}$ and $F^s\widehat{\Omega}^k_{b,0}\to
F^s\widehat{\Omega}^{k+1}_{b,0}$ defined in section 4 give complexes
$\widetilde{\Omega}^\bullet_{b,0}$ and $\widehat{\Omega}^\bullet_{b,0}$
As in (4.6) and (4.7) we have
\begin{equation}
H^0(\widehat{\Omega}^\bullet_{b,0}) = \tilde{\Lambda}_0
\end{equation}
and
\begin{equation}
H^k(\widehat{\Omega}^\bullet_{b,0}) =
H^{k-1}(\widetilde{\Omega}^\bullet_{b,0}) \quad\text{for $k\geq 1$,}
\end{equation}
and as in (4.8) we have a short exact sequence of complexes
\begin{equation}
0\to \widetilde{\Omega}^\bullet_{b,0}\to \Omega^\bullet_{b,0}\to
\widehat{\Omega}^\bullet_{b,0}\to 0.
\end{equation} 

It is straightforward
to check from (2.14) and (4.1) that 
\begin{equation}
\theta\circ q\alpha_k = \alpha_{k-1}\circ\theta
\end{equation}
for $k=1,\dots,n+r+1$.  It follows that $\alpha_k$ is stable on
$\widetilde{\Omega}^k_{b,0}$, hence the Frobenius structure on
$\Omega^\bullet_{b,0}$ induces a Frobenius structure on
$\widetilde{\Omega}^\bullet_{b,0}$ and the inclusion
$\widetilde{\Omega}^\bullet_{b,0} \hookrightarrow
\Omega^\bullet_{b,0}$ is a morphism of complexes with Frobenius
structure.  For the Frobenius structure on
$\widehat{\Omega}^\bullet_{b,0}$, we define
$\hat{\alpha}_k:\widehat{\Omega}^k_{b,0}\to \widehat{\Omega}^k_{b,0}$
as follows.  For $k=0$ we have $\widehat{\Omega}^0_{b,0} =
\tilde{\Lambda}_0$, and we define $\hat{\alpha}_0$ to be
multiplication by $q^{n+r+1}$.  For $k>0$ we have
$\widehat{\Omega}^k_{b,0} = \widetilde{\Omega}^{k-1}_{b,0}$ and we
define $\hat{\alpha}_k$ to be $q^{-1}\alpha_{k-1}$.  Then (6.5)
implies that the map
$\Omega^\bullet_{b,0}\to\widehat{\Omega}^\bullet_{b,0}$ is a morphism
of complexes with Frobenius structure.  From (6.4) we then get an
exact sequence of cohomology spaces with Frobenius structure
\begin{equation}
\cdots\to H^k(\widetilde{\Omega}^\bullet_{b,0})\to
H^k(\Omega^\bullet_{b,0})\to H^k(\widehat{\Omega}^\bullet_{b,0})
\xrightarrow{\delta}
H^{k+1}(\widetilde{\Omega}^\bullet_{b,0})\to\cdots,
\end{equation}
where $\delta$ denotes the connecting homomorphism.  Note in
particular that the identification (6.3) is {\em not} an isomorphism of
Frobenius modules.  One has instead
\begin{equation}
\det(I-qt\hat{\alpha}_k\mid H^k(\widehat{\Omega}^\bullet_{b,0})) = 
\det(I-t\alpha_{k-1}\mid H^{k-1}(\widetilde{\Omega}^\bullet_{b,0})).
\end{equation}

From Propositions 4.12 and 7.7 we get the following result.
\begin{proposition}
Let $0\leq k<2r$.  For $k$ even, $H^k(\widetilde{\Omega}^\bullet_{b,0}) =
H^{k+1}(\widehat{\Omega}^\bullet_{b,0}) = 0$, 
and for $k$ odd, $H^k(\widetilde{\Omega}^\bullet_{b,0}) =
H^{k+1}(\widehat{\Omega}^\bullet_{b,0})$ is a one-dimensional vector
space with basis $[\eta_{(k+1)/2}]$.
\end{proposition}

By Lemma 4.10 we have $\delta([\eta_i]) = [\eta_{i+1}]$, so for $k$
even, $0\leq k\leq 2r-2$, the connecting homomorphism $\delta$ of
(6.6) is an isomorphism of one-dimensional vector spaces with
Frobenius and we get 
\begin{equation}
\det(I-t\hat{\alpha}_k\mid H^k(\widehat{\Omega}^\bullet_{b,0})) =
\det(I-t\alpha_{k+1}\mid H^{k+1}(\widetilde{\Omega}^\bullet_{b,0})).
\end{equation}
Combining this with the observation (6.7) gives for $k$ even, $0\leq
k\leq 2r-4$, 
\begin{equation}
\det(I-t\hat{\alpha}_k\mid H^k(\widehat{\Omega}^\bullet_{b,0})) =
\det(I-qt\hat{\alpha}_{k+2}\mid H^{k+2}(\widehat{\Omega}^\bullet_{b,0})). 
\end{equation}
From the definition of $\hat{\alpha}_0$ we have
$\det(I-t\hat{\alpha}_0\mid H^0(\widehat{\Omega}^\bullet_{b,0})) =
(1-q^{n+r+1}t)$, hence for $k$ even, $0\leq k\leq 2r-2$, we have by (6.10)
that
\begin{equation}
\det(I-t\hat{\alpha}_k\mid H^k(\widehat{\Omega}^\bullet_{b,0})) =
(1-q^{n+r+1-(k/2)}t).
\end{equation}
From (6.9), we then get for these same $k$ that
\begin{equation}
\det(I-t\alpha_{k+1}\mid H^{k+1}(\widetilde{\Omega}^\bullet_{b,0})) =
(1-q^{n+r+1-(k/2)}t).
\end{equation}
Finally, taking $k=2r-2$ in (6.12) and using (6.7) gives
\begin{equation}
\det(I-t\hat{\alpha}_{2r}\mid H^{2r}(\widehat{\Omega}^\bullet_{b,0}))
= (1-q^{n+1}t).
\end{equation}

By Propositions 3.12 and 7.7, the cohomology class $[\xi_r]$ is a basis for
$H^{2r}(\Omega^\bullet_{b,0})$, and by Propositions 4.12 and 7.7, the
cohomology class $[\eta_r]$ is a basis for
$H^{2r}(\widehat{\Omega}^\bullet_{b,0})$.   
\begin{proposition}
Let $r<n$.  Relative to the bases $[\xi_r]$ for
$H^{2r}(\Omega^\bullet_{b,0})$ and $[\eta_r]$ for
$H^{2r}(\widehat{\Omega}^\bullet_{b,0})$, the map
\[ \theta:H^{2r}(\Omega^\bullet_{b,0})\to
H^{2r}(\widehat{\Omega}^\bullet_{b,0}) \]
is multiplication by $(-1)^{r(r-1)/2}d_1\cdots d_r$.  In particular,
this map is an isomorphism of Frobenius modules, hence
\[ \det(I-t\alpha_{2r}\mid H^{2r}(\Omega^\bullet_{b,0})) =
(1-q^{n+1}t). \]
\end{proposition}

{\it Proof}.  The first assertion is an immediate consequence of
Proposition 4.19.  The second is an immediate consequence of the first
and equation (6.13).

From Propositions 4.23 and 7.7 we get the following result.
\begin{proposition}
For $2r\leq k<n+r$, $H^k(\widetilde{\Omega}^\bullet_{b,0}) =
H^{k+1}(\widehat{\Omega}^\bullet_{b,0}) = 0$.
\end{proposition}

We can now state the main consequence of our cohomological computations.
\begin{theorem}
\begin{multline*}
L({\bf A}^{n+r+1},\Psi,F;t) = \\
(1-q^{n+1}t)^{-1}\biggl(\frac{\det(I-t\alpha_{n+r+1}\mid 
H^{n+r+1}(\Omega^\bullet_{b,0}))}{\det(I-qt\alpha_{n+r+1}\mid 
H^{n+r+1}(\Omega^\bullet_{b,0}))}\biggr)^{(-1)^{n+r}}. 
\end{multline*}
\end{theorem}

{\it Proof}.  From (2.20) and Proposition 3.6, if $r<n$ then
\begin{align*}
L({\bf A}^{n+1+r},&\Psi,F;t) = \\ 
\begin{split}
& \det(I-t\alpha_{2r}\mid H^{2r}(\Omega^\bullet_{b,0}))^{-1}
\det(I-t\alpha_{n+r}\mid
H^{n+r}(\Omega^\bullet_{b,0}))^{(-1)^{n+r+1}}\\ 
\cdot&\det(I-t\alpha_{n+r+1}\mid
H^{n+r+1}(\Omega^\bullet_{b,0}))^{(-1)^{n+r}},  
\end{split}
\end{align*}
while if $r=n$ (so that $2r=n+r$), then
\begin{multline*}
L({\bf A}^{n+1+r},\Psi,F;t) = \\
\det(I-t\alpha_{n+r}\mid
H^{n+r}(\Omega^\bullet_{b,0}))^{-1}\det(I-t\alpha_{n+r+1}\mid
H^{n+r+1}(\Omega^\bullet_{b,0})). 
\end{multline*}
If $r<n$, then by Proposition 6.14, Theorem 6.16 reduces to proving
that
\begin{equation}
\det(I-t\alpha_{n+r}\mid
H^{n+r}(\Omega^\bullet_{b,0})) =
\det(I-qt\alpha_{n+r+1}\mid H^{n+r+1}(\Omega^\bullet_{b,0})).
\end{equation}
If $r=n$, then Theorem 6.16 reduces to proving that
\begin{equation}
\det(I-t\alpha_{n+r}\mid
H^{n+r}(\Omega^\bullet_{b,0})) = (1-q^{n+1}t)
\det(I-qt\alpha_{n+r+1}\mid H^{n+r+1}(\Omega^\bullet_{b,0})).
\end{equation}

Note that since $\widetilde{\Omega}^k_{b,0} = 0$ for $k>n+r$, the
exact sequence (6.6) gives 
\begin{equation}
H^{n+r+1}(\Omega^\bullet_{b,0})\cong
H^{n+r+1}(\widehat{\Omega}^\bullet_{b,0})
\end{equation}
Suppose $r<n-1$.  Then
\[ H^{n+r-1}(\widehat{\Omega}^\bullet_{b,0}) =
H^{n+r}(\widehat{\Omega}^\bullet_{b,0})=0 \]
by Proposition 6.15.  Using this in (6.6) gives
\begin{equation}
H^{n+r}(\Omega^\bullet_{b,0})\cong
H^{n+r}(\widetilde{\Omega}^\bullet_{b,0}).
\end{equation}
Equation (6.17) now follows from (6.19), (6.20), and (6.7).  If
$r=n-1$, then $H^{n+r}(\widehat{\Omega}^\bullet_{b,0})=0$ by
Proposition 6.15, so (6.6) gives an exact sequence
\[ H^{2r}(\widehat{\Omega}^\bullet_{b,0})\xrightarrow{\delta}
H^{n+r}(\widetilde{\Omega}^\bullet_{b,0})\to
H^{n+r}(\Omega^\bullet_{b,0})\to 0. \]
By Proposition 6.14, the image of $\delta$ is spanned by
\[ \delta([\theta(\xi_r)]) = [D(\xi_r)] = 0, \]
i.~e., $\delta$ is the zero map, so (6.20) holds in this case also.
Equation (6.17) now follows as in the case $r<n-1$.  Finally, suppose
$r=n$.  Then $H^{n+r-1}(\widehat{\Omega}^\bullet_{b,0}) = 0$ by
Proposition 6.8, so (6.6) gives an exact sequence
\begin{equation}
0\to H^{n+r}(\widetilde{\Omega}^\bullet_{b,0})\to
H^{n+r}(\Omega^\bullet_{b,0})\xrightarrow{\theta}
H^{n+r}(\widehat{\Omega}^\bullet_{b,0})\to 0. 
\end{equation}
By (6.7) and (6.19) we have
\begin{equation}
\det(I-t\alpha_{n+r}\mid H^{n+r}(\widetilde{\Omega}^\bullet_{b,0}))
= \det(I-qt\alpha_{n+r+1}\mid H^{n+r+1}(\Omega^\bullet_{b,0})).
\end{equation}
Equation (6.18) then follows from (6.21), (6.22) and (6.13).  This
completes the proof of the theorem.

\begin{corollary}
Let $P(t)$ be as defined in section $1$.  Then
\[ P(q^rt) = \det(I-t\alpha_{n+r+1}\mid
H^{n+r+1}(\Omega^\bullet_{b,0})). \]
\end{corollary}

{\it Proof}.  Let $g(t)\in 1+ t\tilde{\Lambda}_0[t]$ and define $r(t)
= g(t)/g(qt)\in 1+t\tilde{\Lambda}_0[[t]]$.  Then $\prod_{i=0}^{m-1}
r(q^it) = g(t)/g(q^mt)$.  It follows that
\[ g(t) = \lim_{m\to\infty}\prod_{i=0}^{m-1}r(q^it), \]
in the sense that the coefficients of the power series on the
right-hand side converge term-by-term to the coefficients of $g(t)$.
In particular, $g(t)$ is uniquely determined by $r(t)$.  The assertion
of the corollary now follows from (2.3) and Theorem 6.16.

By Corollary 6.23, Theorem 1.1 is equivalent to the following result.
\begin{theorem}
Suppose that $f_1=\cdots=f_r=0$ defines a smooth complete intersection
$X$ in ${\bf P}^n$.  Then the Newton polygon of
$\det(I-t\alpha_{n+r+1}\mid H^{n+r+1}(\Omega^\bullet_{b,0}))$ with
respect to ${\rm ord}_q$ lies on or above the Newton polygon with
respect to ${\rm ord}_q$ of the polynomial $\prod_{e=r}^n
(1-q^et)^{h_e}$. 
\end{theorem}

We begin with a reduction step.  In (2.15), 
we defined a $\tilde{\Lambda}_1$-linear endomorphism $\beta_{n+r+1}$ of
$H^{n+r+1}(\Omega^\bullet_{b,0})$ such that $\alpha_{n+r+1} =
(\beta_{n+r+1})^a$.  Let ${\rm ord}$ denote the $p$-adic valuation
normalized by ${\rm ord}\;p = 1$.  By \cite[Lemma 7.1]{DW2} we have
the following.
\begin{lemma}
The Newton polygon of $\det(I-t\alpha_{n+r+1}\mid
H^{n+r+1}(\Omega^\bullet_{b,0}))$ with respect to the valuation
${\rm ord}_q$ is obtained from the Newton polygon of
$\det(I-t\beta_{n+r+1}\mid H^{n+r+1}(\Omega^\bullet_{b,0}))$ with
respect to the valuation ${\rm ord}$ by shrinking ordinates and
abscissas by a factor of $1/a$. 
\end{lemma}

Theorem 6.24 is thus equivalent to the following result.
\begin{theorem}
Suppose that $f_1=\cdots=f_r=0$ defines a smooth complete intersection
$X$ in ${\bf P}^n$.  Then the Newton
polygon of $\det(I-t\beta_{n+r+1}\mid
H^{n+r+1}(\Omega^\bullet_{b,0}))$ with respect to ${\rm ord}$ lies
on or above the Newton polygon with respect to ${\rm ord}$ of the
polynomial $\prod_{e=r}^n (1-p^et)^{ah_e}$.
\end{theorem}

{\it Proof of Theorem $6.26$}.  Let $\{\gamma_m\}_{m=1}^a$ be an
integral basis for $\tilde{\Lambda}_0$ over $\tilde{\Lambda}_1$.  By
Theorem 5.18 and the definition of an integral basis, the cohomology
classes 
\begin{equation}
[\gamma_m\xi^{(e)}_l], \quad e=r,\dots,n, \quad l=1,\dots,h_e, \quad
m=1,\dots,a,
\end{equation}
form a basis for $H^{n+r+1}(\Omega^\bullet_{b,0})$ as
$\tilde{\Lambda}_1$-vector space.  It is straightforward to check from
the definitions that
\[ \beta(F^sC(b/p))\subseteq F^sC(b) \]
for all $s\in{\bf Z}$.  Using this, one checks that 
\begin{equation}
\beta_{n+r+1}(\gamma_m\xi^{(e)}_l) \in
F^{Mbe(p-1)/p}\Omega^{n+r+1}_{b,0}.
\end{equation}
By Corollary 7.8,
$[\beta_{n+r+1}(\gamma_m\xi^{(e)}_l)]$ is a $\tilde{\Lambda}_1$-linear 
combination of the $[\gamma_{m'}\xi^{(e')}_{l'})]$ with coefficients
in $\pi^{Mbe(p-1)/p}{\mathcal O}_{\tilde{\Lambda}_1}$.  This says that
in the matrix of $\beta_{n+r+1}$ relative to the basis (6.27), the
column corresponding to $[\gamma_m\xi^{(e)}_l]$ has all entries
divisible by $\pi^{Mbe(p-1)/p}$.  This implies that the Newton polygon
of $\det(I-t\beta_{n+r+1}\mid H^{n+r+1}(\Omega^\bullet_{b,0}))$ 
with respect to the valuation ord lies on or above the Newton
polygon with respect to the valuation ord of the polynomial
\[ \prod_{e=r}^n (1-\pi^{Mbe(p-1)/p}t)^{ah_e}. \]
But $\det(I-t\beta_{n+r+1}\mid H^{n+r+1}(\Omega^\bullet_{b,0}))$ is
independent of $b$ by Corollary 6.23, so we may take the limit as $b\to
p/(p-1)$ to conclude (recall that $\pi^M = p$) that its Newton polygon
lies on or above the Newton polygon of $\prod_{e=r}^n (1-p^et)^{ah_e}$. 

\section{Appendix}

In this section we collect (with references) some basic results on
``lifting'' cohomology from characteristic $p$ to characteristic zero.
Let ${\mathcal O}$ be a complete discrete valuation ring with
uniformizer $\pi$.  Call an ${\mathcal O}$-module $M$ {\it flat\/} if
multiplication by $\pi$ is injective and call $M$ separated if
$\bigcap_{j=1}^\infty \pi^jM = 0$.  A separated ${\mathcal O}$ -module
$M$ has an obvious metric space structure with the
$\{\pi^jM\}_{j=1}^\infty$ forming a fundamental system of neighborhoods
of $0$.  Call $M$ ${\mathcal O}$-{\it complete\/} if it is complete in
this metric.  Let 
\[ C^\bullet = \{ 0\to C^0\xrightarrow{\partial} C^1\to\dots\} \]
be a complex of flat, separated, ${\mathcal O}$-complete ${\mathcal
O}$-modules with ${\mathcal O}$-linear boundary maps.  Let 
\[ \bar{C}^\bullet = \{ 0\to \bar{C}^0\xrightarrow{\bar{\partial}}
\bar{C}^1\to\dots\} \]  
be the complex obtained by reducing $C^\bullet$ modulo $\pi$, i.~e.,
$\bar{C}^i = C^i/\pi C^i$ and the boundary maps of $\bar{C}^\bullet$
are those induced by the boundary maps of $C^\bullet$.

The first assertion of the following result is \cite[Theorem
A.1(a)]{AS0}.  The second assertion follows from the proof of
\cite[Theorem A.1(a)]{AS0}. 
\begin{proposition}
If $H^i(\bar{C}^\bullet) = 0$ for some $i$, then $H^i(C^\bullet)=0$.
More precisely, if $\omega\in C^i$ satisfies $\partial(\omega) = 0$
and if $\bar{\omega} = \bar{\partial}(\eta)$ for some $\eta\in
\bar{C}^{i-1}$ (where $\bar{\omega}$
denotes the image of $\omega$ in $\bar{C}^i$), then there exists
$\xi\in C^{i-1}$ such that $\partial(\xi) = \omega$ and $\bar{\xi} =
\eta$.  
\end{proposition}

The next result is \cite[Lemma 4.1]{AS1}.
\begin{proposition}
If $H^i(\bar{C}^\bullet)$ is of finite dimension $d$ over ${\mathcal
O}/(\pi)$ and multiplication by $\pi$ is injective on $H^i(C^\bullet)$
and $H^{i+1}(C^\bullet)$, then $H^i(C^\bullet)$ is a free ${\mathcal
O}$-module of rank~$d$.  Furthermore, if $\xi_1,\dots,\xi_d\in C^i$
satisfy (a) $\partial(\xi_j) = 0$ for $j=1,\dots,d$ and (b) the
cohomology classes $[\bar{\xi}_1],\dots,[\bar{\xi}_d]$ form a basis
for $H^i(\bar{C}^\bullet)$, then the cohomology classes
$[\xi_1],\dots,[\xi_d]$ form a basis for $H^i(C^\bullet)$.  
\end{proposition}

We shall also need some additional results.
\begin{proposition}
Let $\{\xi_j\}_{j=1}^N\subseteq C^i$ satisfy (a) $\partial(\xi_j)=0$
for all $j$ and (b) the cohomology classes $\{[\bar{\xi}_j]\}_{j=1}^N$
span $H^i(\bar{C}^\bullet)$.  Then the cohomology classes
$\{[\xi_j]\}_{j=1}^N$ span $H^i(C^\bullet)$.
\end{proposition}

{\it Proof}.  Let $\omega\in C^i$ with $\partial(\omega) = 0$.  Then
$\bar{\partial}(\bar{\omega})=0$, so there exist
$\alpha^{(0)}_j\in{\mathcal O}$ and $\eta_0\in C^{i-1}$ such that
\begin{equation}
\omega = \sum_{j=1}^N \alpha^{(0)}_j\xi_j + \partial(\eta_0) +
\pi\omega_0
\end{equation}
for some $\omega_0\in C^i$.  Suppose that for some $m\geq 0$ we have
\begin{equation}
\omega = \sum_{j=1}^N \alpha^{(m)}_j\xi_j + \partial(\eta_m) +
\pi^{m+1}\omega_m
\end{equation}
for some $\alpha^{(m)}_j\in{\mathcal O}$, $\eta_m\in C^{i-1}$, and
$\omega_m\in C^i$ with
\[ \alpha^{(m)}_j-\alpha^{(m-1)}_j\in \pi^m{\mathcal O} \quad
\text{and} \quad \eta_m-\eta_{m-1}\in \pi^m C^{i-1}. \]
Equation (7.5) implies that $\partial(\omega_m) = 0$, so as in (7.4)
we have
\begin{equation}
\omega_m = \sum_{j=1}^N \alpha'_j\xi_j + \partial(\eta') +
\pi\omega_{m+1}
\end{equation}
for some $\alpha'_j\in{\mathcal O}$, $\eta'\in C^{i-1}$, and
$\omega_{m+1}\in C^i$.  Put
\[ \alpha_j^{(m+1)} = \alpha_j^{(m)} + \pi^{m+1}\alpha'_j \quad
\text{and} \quad \eta_{m+1} = \eta_m + \pi^{m+1}\eta'. \]
Substituting (7.6) into (7.5) gives
\[ \omega = \sum_{j=1}^N \alpha^{(m+1)}_j\xi_j + \partial(\eta_{m+1}) +
\pi^{m+2}\omega_{m+1} \]
with 
\[ \alpha_j^{(m+1)}-\alpha_j^{(m)}\in \pi^{m+1}{\mathcal O}
\quad\text{and}\quad \eta_{m+1}-\eta_m\in \pi^{m+1}C^{i-1}. \]
It follows that each sequence $\{\alpha_j^{(m)}\}_{m=0}^\infty$
converges to an element $\alpha_j\in{\mathcal O}$ and
$\{\eta_j\}_{j=0}^\infty$ converges to an element $\eta\in C^{i-1}$
satisfying
\[ \omega = \sum_{j=1}^N \alpha_j\xi_j + \partial(\eta), \]
which shows that the $[\xi_j]$ span $H^i(C^\bullet)$.

Let $\Lambda$ be the quotient field of ${\mathcal O}$ and suppose that
$D^\bullet = \{0\to D^0\xrightarrow{\partial} D^1\to 
\dots\}$ is a complex of $\Lambda$-vector spaces containing
$C^\bullet$ as a subcomplex such that
\[ D^\bullet = \bigcup_{s\in{\bf Z}} \pi^sC^\bullet. \]
\begin{proposition}
Suppose that $H^i(C^\bullet) = H\oplus H'$, where $H$ is a free
${\mathcal O}$-module with basis the cohomology classes
$[\xi_1],\dots,[\xi_d]$ and $H'$ is a torsion ${\mathcal O}$-module
spanned by the cohomology classes $[\eta_1],\dots,[\eta_e]$.  Then the
cohomology classes $[\xi_1],\dots,[\xi_d]$ form a basis for
$H^i(D^\bullet)$ as $\Lambda$-vector space.
\end{proposition}

{\it Proof}.  Since multiplication by $\pi^s$ is an isomorphism from
$C^\bullet$ to $\pi^sC^\bullet$, it follows that $H^i(\pi^sC^\bullet)
= H_s\oplus H'_s$, where $H_s$ is a free ${\mathcal O}$-module with
basis the cohomology classes $[\pi^s\xi_1],\dots,[\pi^s\xi_d]$ and
$H'_s$ is a torsion ${\mathcal O}$-module spanned by the cohomology
classes $[\pi^s\eta_1],\dots,[\pi^s\eta_e]$.  Let $\xi\in D^i$ with
$\partial(\xi) = 0$.  Since $\xi\in\pi^sC^i$ for some $s$ we have
\[ \xi = \sum_{k=1}^d c_k\pi^s\xi_k + \sum_{l=1}^e c'_l\pi^s\eta_l +
\partial(\zeta) \]
for some $c_k,c'_l\in{\mathcal O}$ and some $\zeta\in\pi^sC^{i-1}$.
Since the $[\pi^s\eta_l]$ are torsion elements, it follows that there
exists a positive integer $t$ such that
\[ \pi^t\pi^s\eta_l = \partial(\zeta_l) \]
for $l=1,\dots,e$, where $\zeta_l\in \pi^sC^{i-1}$.  Substitution then
gives
\[ \xi = \sum_{k=1}^d c_k\pi^s\xi_k + \partial(\zeta + \sum_{l=1}^e
c'_l\pi^{-t}\zeta_l), \]
which shows that the $[\xi_k]$ span $H^i(D^\bullet)$.   

Suppose there are $c_k\in\Lambda$ and $\zeta\in D^{i-1}$ such that
\[ \sum_{k=1}^d c_k\xi_k = \partial(\zeta). \]
There exists an integer $s$ such that $\pi^sc_k\in{\mathcal O}$ for
all $k$ and $\zeta\in\pi^{-s}C^{i-1}$.  Thus we have
\[ \sum_{k=1}^d (\pi^sc_k)(\pi^{-s}\xi_k) = \partial(\zeta). \]
But if some $c_k\neq 0$, this contradicts the fact that the
$[\pi^{-s}\xi_k]$ are a basis for the free ${\mathcal O}$-module
$H_{-s}$.  This proves that the $[\xi_k]$ are linearly independent in
$H^i(D^\bullet)$.  

The proof of Proposition 7.7 shows that the following assertion holds.
\begin{corollary}
Let $\xi\in D^i$ with $\partial(\xi) = 0$.  If $\xi\in\pi^s C^i$, then
\[ [\xi] = \sum_{k=1}^d c_k[\xi_k] \]
with $c_k\in\pi^s{\mathcal O}$ for all $k$.
\end{corollary}

\end{document}